\documentclass[12pt]{amsart}

\usepackage[matrix,arrow,curve,frame]{xy}
\usepackage{amsmath,amsthm,amssymb,enumerate}
\usepackage{latexsym}
\usepackage{amscd}
\usepackage[colorlinks=true]{hyperref}
\usepackage{euscript}

\newtheorem{thm}[subsection]{Theorem}
\newtheorem{defn}[subsection]{Definition}

\newtheorem{claim}[subsection]{Claim}
\newtheorem{corr}[subsection]{Corollary}
\newtheorem{lemma}[subsection]{Lemma}

\newtheorem{remark}[subsection]{Remark}

\theoremstyle{definition}

\newcommand{\cat}{\mathcal}

\newcommand{\R}{\mathbb R}
\newcommand{\Q}{\mathbb Q}
\newcommand{\Z}{\mathbb Z}
\newcommand{\F}{\mathbb F}

\newcommand{\C}{\mathbb C}

\newcommand{\T}{\mathbb T}

\newcommand{\gG}{\mbox{\german g}}

\newcommand{\gB}{\mathfrak{b}}

\newcommand{\gC}{{\cat S}}

\newcommand{\Oo}{{\cat O}}

\DeclareMathOperator{\colim}{colim}
\DeclareMathOperator{\hocolim}{hocolim}

\DeclareMathOperator{\Hom}{Hom}

\newfont{\german}{eufm10}

 \DeclareMathOperator{\Img}{Im}

 \DeclareMathOperator{\Ker}{Ker}

 \DeclareMathOperator{\BK}{BK}
\DeclareMathOperator{\BN(T)}{BN(T)}
\DeclareMathOperator{\BNJ}{BN_J(T)}
\DeclareMathOperator{\BW}{BW}

\DeclareMathOperator{\BT}{BT}

\DeclareMathOperator{\BG}{BG}
\DeclareMathOperator{\GrN}{G/N}
\DeclareMathOperator{\KNT}{K/N(T)}
\DeclareMathOperator{\KNJT}{K_J/N_J(T)}
\DeclareMathOperator{\KNJKT}{K_{J_k}/N_{J_k}(T)}
\DeclareMathOperator{\KNJIT}{K_{J_i}/N_{J_i}(T)}
\DeclareMathOperator{\KNST}{K_S/N_S(T)}
\DeclareMathOperator{\To}{T}
\DeclareMathOperator{\W}{W}
\DeclareMathOperator{\K}{K}
\DeclareMathOperator{\KT}{K/T}
\DeclareMathOperator{\BS}{BS^1}
\DeclareMathOperator{\So}{S^1}
\DeclareMathOperator{\No}{N}
\DeclareMathOperator{\B}{B}
\DeclareMathOperator{\Hr}{H}
\DeclareMathOperator{\ET}{ET}
\DeclareMathOperator{\Ro}{R}
\DeclareMathOperator{\X}{X}
\DeclareMathOperator{\Gr}{G}
\DeclareMathOperator{\Y}{Y}
\DeclareMathOperator{\Zo}{Z}
 \DeclareMathOperator{\NoT}{N(T)}
 \DeclareMathOperator{\Tor}{Tor}
 \DeclareMathOperator{\KIT}{K_i/T}
 \DeclareMathOperator{\KJT}{K_j/T}
 \DeclareMathOperator{\KKI}{K/K_i}
 \DeclareMathOperator{\KKJ}{K/K_j}
 \DeclareMathOperator{\TT}{\mathcal{T}}

 \newcommand{\gH}{\mathfrak{h}}

 \newcommand{\tens} {\otimes}
 \newcommand{\iso} {\cong}

\begin{document}
\pagestyle{plain}

\title
{On the Topology of Kac-Moody groups}
\author{Nitu Kitchloo}
\address{Department of Mathematics, Johns Hopkins University, Baltimore, USA}
\email{nitu@math.jhu.edu}
\thanks{Nitu Kitchloo is supported in part by NSF through grant DMS
  1005391.}

\date{\today}


{\abstract

\noindent
We study the topology of spaces related to Kac-Moody groups. Given a split Kac-Moody group over $\C$, let $\K$ denote the unitary form with maximal torus $\To$ having normalizer $\No(\To)$. In this article we study the cohomology of the flag manifold $\KT$ as a module over the Nil-Hecke algebra, as well as the (co)homology of $\K$ as a Hopf algebra. In particular, if $\F$ has positive characteristic, we show that $\Hr_*(\K,\F)$ is a finitely generated algebra, and that $\Hr^*(\K,\F)$ is finitely generated only if $\K$ is a compact Lie group . We also study the stable homotopy type of the classifying space $\BK$ and show that it is a retract of the classifying space $\BN(T)$ of $\No(\To)$. We illustrate our results with the example of rank two Kac-Moody groups.}
\maketitle

\tableofcontents
\section{Introduction}

\noindent
In this paper, we deal with a class of topological groups known as Kac-Moody groups \cite{KP2}. By a Kac-Moody group, we shall mean the unitary form of a split Kac-Moody group over $\C$. We refer the reader to \cite{Ku} for a beautiful treatment of the subject. These groups form a natural extension of the class of compact Lie groups, and share many of their properties. They are known to contain the class of (polynomial) loop groups, which go by the name of affine Kac-Moody groups. With the exception of compact Lie groups, Kac-Moody groups over $\C$ are not even locally compact. Hence geometric arguments used to study the topology of compact Lie groups and their flag varieties no longer extend to general Kac-Moody groups. This led Kac-Peterson to construct of a whole new set of techniques applicable in this context. Underlying these techniques is a collection of (annhilation) operators that generate a ring $\Oo$ \cite{K1,KK} given by a deformation of the group ring of the Weyl group. The cohomology rings of the flag varieties admit an action of $\Oo$ which intertwines the action of the Weyl group on the flag varieties. This induces a very rich structure that can be exploited to prove various structure theorems about flag varieties \cite{Ki,K1,Ku}, as well as the Kac-Moody group itself \cite{K1}. Further techniques were developed by the author in \cite{Ki} in order to study the classifying spaces of these groups. In \cite{Ki, Ki2} it was shown that these classifying spaces can be decomposed as a suitable homotopy colimit of classifying spaces of compact Lie groups. 

\medskip
\noindent
The main purpose of this paper is to prove various new results regarding the topology of Kac-Moody groups. The first set of new results pertain to the (co)homology rings of Kac-Moody groups and their flag varieties with coefficients in a field of positive characteristic.  Next, we study the classifying space of a Kac-Moody group and prove various new results. In particular,  we show that the stable homotopy type of the classifying space is a retract of the classifying space of the normalizer of the maximal torus.

\medskip
\noindent
In addition to the results mentioned above, we will also provide independent proofs of various homological structure theorems, many of which are known to the experts, but whose proofs are either absent in print or exist in weaker generality. Some of these results were stated in \cite{K1}, but the proofs never made it to print. 

\subsection{Acknowledgements:}

\noindent

\noindent
Several of these results presented here are borrowed from the author's thesis \cite{Ki}. The author would like to thank his advisor Haynes Miller for introducing him to Kac-Moody groups. It is a pleasure to thank Victor Kac and Dale Peterson for all the knowledge they have shared with him. In particular, the author benefited greatly from his interactions with Dale Peterson at M.I.T. in the spring of 1997. Also, the author thanks Kasper Andersen for detailed comments on this document. And finally, we would like to thank the referee for carefully reading this document, and helping us clean up the arguments.

\section{Background and Statement of New Results}

\noindent
Kac-Moody groups have been extensively studied and much is known about their general structure, representation theory and topology \cite{K1,KP2,KP3,Ki,Ku,KW,T} (see \cite{Ku} for a modern prespective). Their construction begins with a finite integral matrix $A = (a_{ij})_{i,j \in I}$ with the properties that $a_{ii} = 2$ and $a_{ij} \leq 0$ for $i\neq j$. Moreover, we demand that $a_{ij} = 0$ if and only if $a_{ji} = 0$. These conditions define a {\em Generalized Cartan Matrix}. 

\medskip
\noindent
Given a generalized Cartan matrix $A$, one may construct a complex Lie algebra $\gG (A)$ using the Harishchandra-Serre relations. This Lie algebra contains a finite dimensional Cartan subalgebra $\gH$ that admits an integral form $\gH_{\Z}$ and a real form $\gH_{\R} = \gH_{\Z}\otimes \R$. The lattice $\gH_{\Z}$ contains a finite set of primitive elements $h_i, i \in I$ called "simple coroots".  Similarly, the dual lattice $\gH_{\Z}^*$ contains a special set of elements called "simple roots" $\alpha_i, i \in I$. One may decompose $\gG(A)$ under the adjoint action of $\gH$ to obtain a triangular form as in the classical theory of semisimple Lie algebras. Let $\eta_{\pm}$ denote the positive and negative "nilpotent" subalgebras respectively, and let $\gB_{\pm} = \gH \oplus \eta_{\pm}$ denote the corresponding "Borel" subalgebas. The structure theory of $\gG(A)$ leads to a construction (in much the same way that Chevalley groups are constructed), of a topological group $\Gr(A)$ called the (minimal, split) Kac-Moody group over the complex numbers. The group $\Gr(A)$ supports a canonical anti-linear involution $\omega$, and one defines the unitary form $\K(A)$ as the fixed group $\Gr(A)^{\omega}$. The inclusion map $\K(A) \subset \Gr(A)$ is a homotopy equivalence. It is more convenient for your purposes to work with the group $\K(A)$ in this article. 

\medskip
\noindent
Given a subset $J \subseteq I$, one may define a parabolic subalgebra $\gG_J(A) \subseteq \gG(A)$ generated by $\gB_+$ and the root spaces corresponding to the set $J$. For example, $\gG_{\varnothing} (A)= \gB_+$. One may exponentiate these subalgebras to parabolic subgroups $\Gr_J(A) \subset \Gr(A)$. We then define the unitary Levi factors $\K_J(A)$ to be the groups $\K(A) \cap \Gr_J(A)$. Hence $\K_{\varnothing} (A) = \To$ is a torus of rank $2|I| - rk(A)$, called the maximal torus of $\K(A)$. The normalizer $\No(\To)$ of $\To$ in $\K(A)$, is an extension of a discrete group $\W(A)$ by $\To$. The Weyl group $\W(A)$ has the structure of a crystallographic Coxeter group generated by reflections $r_i, i \in I$. $\W(A)$ has a Coxeter presentation given as follows:
\[ \W(A) = \langle r_i, i \in I \, | \, r_i^2=1, \, (r_i \, r_j)^{m_{ij}} = 1 \rangle,  \]
where $m_{ij}$ depends on the generalized Cartan matrix $A = (a_{ij})$, and $m_{ij} =2,3,4,6,\infty$ if $a_{ij}\, a_{ji} = 0,1,2,3,\geq 4$ respectively. For $J \subseteq I$, let $\W_J(A)$ denote the subgroup generated by the corresponding reflections $r_j, j \in J$. The group $\W_J(A)$ is a crystallographic Coxeter group in its own right that can be identified with the Weyl group of $\K_J(A)$. 

\medskip
\noindent
Given a generalized Cartan matrix $A = (a_{ij})_{i,j \in I}$, define a category $\gC(A)$ to be the poset category (under inclusion) of subsets $J \subseteq I$ such that $\K_J(A)$ is a compact Lie group. This is equivalent to demanding that $\W_J(A)$ is a finite group. Notice that $\gC(A)$ contains all subsets of $I$ of cardinality less than two. In particular, $\gC(A)$ is nonempty and has an initial object given by the empty set. The category $\gC(A)$ is also known as the poset of spherical subsets. 
The topology on the group $\K(A)$ is the strong topology generated by the compact subgroups $\K_J(A)$ for $J \in \gC(A)$ (See Appendix). More precisely, $\K(A)$ is the amalgamated product of the compact Lie groups $\K_J(A)$, in the category of topological groups. For an arbitrary subset $L \subseteq I$, the topology induced on homogeneous space of the form $\K(A)/\K_L(A)$ makes it into a CW-complex, with only even cells, indexed by the set of cosets $\W(A)/\W_L(A)$ (see \cite{Ku}). More precisely, the cells can be expressed in terms of the Bruhat decomposition: 
\[ \K(A)/\K_J(A) = \coprod \limits_{w \in \W(A)/\W_J(A)} \B w\,\Gr_J/\Gr_J = \coprod \limits_{w \in \W^J(A)} \B w\B/\B, 
\]
where $\B = \Gr_{\varnothing}(A)$ is the Borel subgroup, $\Gr_J(A)$ are the parabolic subgroups, and $\W^J(A)$ denotes the set of minimal coset representatives for the cosets in $\W(A)/W_J(A)$ (see \cite{H}, Chapter 1.10). The space $\B w\B/\B$ is a subspace homeomorphic to $\C^{l(w)}$, where $l(w)$ denotes the word length of the element $w$ in terms of the generators $r_i$. Hence the cohomology $\Hr^*(\K(A)/\K_J(A))$ has a basis $\{ \delta^w\}_{w \in \W^J(A)}$ where $\delta^w \in \Hr^{2l(w)}(\K(A)/\K_J(A))$ represents the linear dual of the affine cell $\delta_w = [\B w\B/\B] \in \Hr_{2l(w)}(\K(A)/\K_J(A))$. The set $\{\delta^w\}$ is called the Schubert basis. If we are given an inclusion $J < L$, then the maps $\K(A)/\K_J(A) \rightarrow \K(A)/\K_L(A)$ are cellular. Hence we may use the same notation $\delta^w$ for the basis of $\Hr^\ast(\K(A)/\K_J(A))$ regardless of $J$.

\medskip
\noindent
Let $\K(A)/\To \stackrel{\psi} \rightarrow \BT$ be the map that classifies the principal $\To$-bundle $\K(A) \stackrel{\pi} \rightarrow \K(A)/\To$. Let $\F$ be a field and let $J \subset \Hr^*(\BT,\F)$ denote the ideal given by the kernel of $\psi^*$. It has been shown in \cite{K} that $J$ is generated by a regular sequence $\langle \sigma_1, \ldots, \sigma_r \rangle$, with $r \leq rank(\To)$. The ideal $J$ is called the ideal of {\em Generalized Invariants}. Let $S \subseteq \Hr^*(\K(A)/\To,\F)$ denote the subring given by the image of $\psi^*$. It is also shown in \cite{K} that $\Hr^*(\K(A)/\To,\F)$ is a free $S$-module. 

\medskip
\noindent
Let us define a coproduct (introduced by D. Peterson) on the cohomology of $\K(A)/\To$:
\[ \Delta(\delta^w) = \sum_{uv=w} \delta^u \otimes \delta^v, \]
where the sum is taken over all reduced expressions for $w$ i.e. expressions where the minimal word length of $w$ with respect to the generators $r_i$ equals the sum of the minimal word lengths of $u$ and $v$.

\medskip
\noindent
We will prove the following theorem:

\begin{thm} The image in cohomology of the projection map $\K(A) \stackrel{\pi} \rightarrow \K(A)/\To$ is isomorphic to $\Hr^*(\K(A)/\To,\F) \otimes_S \F$. Moreover, this image is a sub Hopf algebra with the coalgebra structure induced via the coproduct $\Delta$ defined above.
\end{thm}

\begin{remark}
The statement of the above theorem was communicated to the author by D. Peterson. The first part of this theorem has been proved in \cite{K}. The second part is stated there without proof. In \cite{K}, V. Kac shows that $\Hr^*(\K(A),\F)$ is free over the image of $\pi^*$, and there is a short exact sequence of algebras:
\[ 1 \rightarrow  \Hr^*(\K(A)/\To,\F)\otimes_S \F \stackrel{\pi^*} \longrightarrow \Hr^*(\K(A),\F) \longrightarrow \Lambda(x_1,\ldots x_r) \rightarrow 1. \]
Here $\Lambda(x_1,\ldots, x_r)$ denotes an exterior algebra on classes $x_i$ of degree given by $deg(\sigma_i)-1$. The number of exterior generators is bounded by the rank of $\T$, and equals it if $\F$ has positive characteristic. In light of the previous theorem, this extension is actually an extension of Hopf algebras. Note that the elements $x_i$ can be choosen to be primitive by \cite{MM} (Thm 7.20). 
\end{remark}
\noindent
Next, will study the homology ring $\Hr_*(\K(A),\F)$ and prove the following:
\begin{thm}
Assume that $\F$ has positive characteristic. Then the Pontrjagin ring $\Hr_*(\K(A),\F)$ is a finitely generated algebra. In addition, $\Hr^*(\K(A),\F)$ and $\Hr^*(\K(A)/\To,\F)$ are finitely generated if and only if $\K(A)$ is a compact Lie group.
\end{thm}

\begin{remark} The Pontrjagin ring $\Hr_*(\K(A),\F)$ will in general be highly non commutative. The structure of the rational Pontrjagin ring: $\Hr_*(\K(A),\Q)$ for a general Kac-Moody group remains unclear to the author. 
\end{remark}

\medskip
\noindent
Now let us consider the classifying space of $\K(A)$, denoted by $\BK(A)$. The study of this space was begun in \cite{Ki} and continued in \cite{BK} and \cite{Ki2}. It was shown in \cite{Ki,BK,Ki2} that $\BK(A)$ can be described as a homotopy colimit in terms of the classifying space $\BK_J(A)$, as $J$ varied through the poset $\gC(A)$. A similar result was proved for the spaces $\BN(T)$ and $\BW(A)$. In this paper, we will use this homotopy decomposition to prove the following:

\begin{thm} Let ${\bf BN(T)}_+$ and ${\bf BK(A)}_+$ denote the suspension spectra of the spaces $\BN(T)_+$ and $\BK(A)_+$ respectively (each endowed with a disjoint base point). Then the canonical map $\BN(T)_+ \rightarrow \BK(A)_+$ admits a stable retraction ${\bf BK(A)}_+ \rightarrow {\bf BN(T)}_+$. In addition, this retraction is equivariant with respect to the center of $\K(A)$. 
\end{thm}

\noindent
And finally, in the appendix we will study the topology on $\K(A)$. In particular, we will show that $\K(A)$ has the homotopy type of a CW complex, and that it is the colimit of the compact Lie groups $\K_J(A)$ in the category of topological groups. 

\section{Cohomology and Integration along the Fiber}
\noindent
In the following sections we supress reference to the generalized Cartan matrix and let $\K$ denote a Kac-Moody group with compact Levi factors $\K_J$, and maximal torus $\To$. Let $\W$ denote the Weyl group. In this section, we will work with cohomology with coefficients in some arbitrary ring $\Ro$.

\medskip
\noindent
Recall the coroots $h_i \in \gH_\Z$. Let $h_i^* \in \gH_{\Z}^*$ denote the dual characters with the property $h_i^*(h_j) = \delta_{ij}$, where $\delta_{ij}$ denotes the Kronecker
symbol. The elements $h_i^*$ need not be unique, but we fix a choice throughout. We shall use the same notation to denote the Euler class $h_i^* \in \Hr^2(\BT)$ of the line bundle $\ET \times_{\To}^{h^\ast_i}\C$ over $\BT$, where the action of $\To$ on $\C$ is obtained by exponentiating $h_i^*$. This induces a map from $\To$ to the product of circles indexed by $i$. By construction, this map is an isomorphism. It follows that the set $\{ h_i^*\}$ is an $\Ro$-basis of $\Hr^2(\BT)$ for any ring $\Ro$. 

\medskip
\noindent
We proceed in a similar fashion with the roots $\alpha_i \in \gH_\Z^*$. We use the same notation to denote the Euler class $\alpha_i \in \Hr^2(\KT)$ of the line bundles $\K \times_{\To}^{\alpha_i} \C$, where $\To$ acts on $\C$ by exponentiating the root $\alpha_i$. Similarly, let $\alpha_i \in \Hr^2(\BT)$ denote the Euler class of the line bundle $\ET \times_{\T}^{\alpha_i}\C$ under the same action. 
 
\noindent
Recall the map $\KT \stackrel{\psi} \rightarrow \BT$ classifying the principal bundle $\K \stackrel{\pi}\rightarrow \KT$. The homomorphism $\psi^*: \Hr^*(\BT) \rightarrow \Hr^*(\KT)$ is called the {\em Characteristic Homomorphism}.
\begin{claim} \label{Euler}
$\psi^*: \Hr^2(\BT) \rightarrow \Hr^2(\KT)$ has the property: $\psi^*(h_i^*) = \delta^{r_i}$, where $\delta^{r_i}$ denotes elements of the Schubert basis of $\Hr^*(\KT)$ corresponding to the generating set $r_i \in \W$. In addition, we have the eqality $\alpha_i = \sum \limits_{j\in I} a_{ji}\, \delta^{r_j}$, where $a_{ji}$ denote the entries of the generalized Cartan matrix.  
\end{claim}
\begin{proof}
The proof follows from the CW decomposition of $\KT$. By by the construction of Kac-Moody groups, there exist injective homomorphisms: $\varphi_i: \operatorname{SU}_2 \rightarrow \K_i \subset \K$ extending the subgroups $\So \subset \To$ obtained by exponentiating the coroot $h_i$. The induced map on the level of flag varieties: $\operatorname{SU_2/S^1} \rightarrow \KT$ is simply the inclusion of the cell corresponding to the Weyl element $r_i$. Now the map $\psi \circ \varphi_i:
\operatorname{SU_2/S^1} \rightarrow \BT$ factors through $\BS$, where $\So \subset \To$ is the circle generated by the coroot $h_i$. It follows that the restriction of the line bundle constructed above: $\ET \times_{\To}^{h^\ast_j}\C$ to $\operatorname{SU_2/S^1}$ via $\psi \circ \varphi_i$ is null if $i \neq j$, and is the Hopf bundle if $i=j$. Now, the map $\psi : \KT \rightarrow \BT$ is two connected, with a basis of $\Hr^2(\KT, \Z)$ given by the classes $\delta^{r_i}$ dual to the Schubert cells $\varphi_i: \operatorname{SU_2/S^1} \subset \KT$.
Restricting along each individual $\varphi_i$ shows that: $\psi^* (h_i^*) = \delta^{r_i}$. Now under the canonical pairing between $\gH_\Z^*$ and $\gH_\Z$, we have the equality: $\alpha_i(h_j) = a_{ji}$ so we get an expression for $\alpha_i$ using the previous part of the claim. 
\end{proof}

\medskip
\noindent
We now proceed to construct certain operators $A_i, i \in I$ acting on $\Hr^*(\KT)$ and $\Hr^*(\BT)$ that were mentioned in the introduction. Let us first introduce the framework:
\begin{defn}
Let $\Hr^*$ denote cohomology with coefficients in a ring $\Ro$ which will be fixed throughout this section. An oriented fibration is a triple $(\pi,n,\tau)$ where:

\medskip
\begin{enumerate}
\item[(1)] $F \rightarrow E \stackrel{\pi}{\rightarrow} B $ is a Serre fibration
\item[(2)] $\Hr^i(F) =0 $ for $i >n$ and $\pi_1(B)$ acts trivially on $\Hr^n(F)$.
\item[(3)] $\tau: \Hr^n(F) \rightarrow \Ro$ is a homomorphism of $\Ro$-modules.
\end{enumerate}

\end{defn}

\medskip
\noindent
Given an oriented fibration $(\pi,n,\tau)$. we can define a homomorphism of $\Ro$-modules:
\[
\int \limits_\tau : \Hr^*(E) \longrightarrow \Hr^{*-n}(B)
\]
called {\em Integration along the fiber} as follows: Consider the Serre spectral sequence for the fibration $F \rightarrow E \stackrel{\pi}{\rightarrow} B$. Then $\int \limits_\tau$ is defined as the composite:
\[
\Hr^*(E) \twoheadrightarrow E_\infty^{*-n,n} \longrightarrow E_2^{*-n,n} = \Hr^{*-n} (B; H^n(F)) 
\stackrel{\tau^*}{\longrightarrow} \Hr^{*-n} (B).
\]
Let us prove an easy claim about this homomorphism:
\medskip
\begin{claim} \label{integration}
Given an oriented fibration $(\pi,n,\tau)$, then $\int_\tau: \Hr^*(E) \rightarrow \Hr^{*-n}(B)$ is a map of $\Hr^*(B)$-modules where $\Hr^*(E)$ is an $\Hr^*(B)$-module via $\pi^*$. In addition, $\int_\tau$ is natural with respect to pullbacks of fibrations. Furthermore, we have $\int \limits_\tau \circ \, \pi^* =0$ if $n>0$.
\end{claim}
\begin{proof}
The first part is obvious using the multiplicative structure of the Serre spectral sequence. To show that $\int \limits_\tau \circ \, \pi^* =0$ if $n>0$, notice first that the image of $\pi^*$ is detected on $E_{\infty}^{*,0}$ via the edge homomorphism. Hence its projection to $E_{\infty}^{*,n}$ is trivial. 
\end{proof}

\begin{defn}
Given $i \in I$, consider the fiber bundle: $\KIT \rightarrow \KT \stackrel{\pi_i}{\rightarrow} \KKI$. Note that $\KIT$ may be canonically identified with $\C P^1$ via the homeomorphism $\varphi_i : \operatorname{SU_2/S^1} \rightarrow \KIT$. Let $\tau_i \in \Hom (\Hr^2(\KIT); \Ro) = \Hom (\Hr^2(\C P^1);\Ro) =\Hr_2(\C P^1)$ be the fundamental class. Note that $(\pi_i,2,\tau_i)$ is an oriented fibration.
Define operators:
\[ A_i: \Hr^*(\KT) \longrightarrow \Hr^{*-2} (\KT) \quad \quad A_i = \pi_i^* \circ \int \limits_{\tau_i}. \]
Similarly for the bundle $\KIT \rightarrow \BT \stackrel{\theta_i}{\rightarrow} \BK_i$, 
we define operators by the same name:
\[ A_i: \Hr^*(\BT) \longrightarrow \Hr^{*-2}(\BT) \quad \quad A_i = \theta_i^* \circ \int \limits_{\tau_i}, \] where $\theta_i$ is the inclusion of the maximal torus. 
\end{defn}

\begin{remark}
Note that the operators $A_i$ are natural with respect to ring homomorphisms $\Ro \rightarrow \Ro'$.
\end{remark}
 \noindent
Now we have a commutative diagram:
\[
\xymatrix{
\KIT        \ar[d] \ar[r]^= & \KIT \ar[d] \\
\KT         \ar[d] \ar[r]^{\psi} & \BT \ar[d] \\
\KKI        \ar[r] & \BK_i
}
\]
From claim \ref{integration} we observe that the operators $A_i$ intertwine with $\psi^\ast$: 
\[
\psi^* \circ A_i = A_i \circ \psi^*.
\]

\section{The Schubert Basis and the Nil Hecke Ring}
\noindent
Recall there is a  CW-decomposition of $\K/\K_J$, for $J \subseteq I$:
\[
\K/\K_J = \coprod \limits_{w \in \W^J} \B w \B/\B. 
\]
We define the Schubert subvarieties of $\K/\K_J$ as closures of cells:
\[
\X^J_w = \overline{\B w\, \B/\B} = \coprod \limits_{w' \in \W^J, w' \leq w} \B w'\B/\B,
\]
where we define a partial order on $\W^J$ by restricting the Bruthat order (\cite{B}, Chapter 2). 

\medskip
\noindent
Notice that $\Hr_{2 l(w)}(\X_w^J)$ is a free $\Ro$-module generated by the largest cell $\delta_w = [\B w\B/\B]$. We define this as the ``fundamental class" $[\X^J_w]$, bearing in mind that this is just an algebraic statement as the spaces $\X^J_w$ are not manifolds by any means. Let $i_w: \X_w^J\rightarrow \K/\K_J$ be the inclusion. 

\medskip
\noindent
\begin{lemma} \label{Schubert} The action of the operators $A_j$ on the Schubert basis is given by:
\[
A_j(\delta^w) = \begin{cases} \delta^{wr_j} \quad \mbox{if } l(wr_j) < l(w)\\
0 \quad \mbox{otherwise}
\end{cases}
\]
\end{lemma}
\begin{proof}
Assume $l(w) = n+1$. By the definition of $A_j$, we known that $A_j(\delta^w)$ is in the image of $\Hr^\ast(\KKJ)$. In particular, we may express $A_j(\delta^w)$ as:
\[ A_j(\delta^w) = \sum_{v \in \W^j, \, l(v) = n} k_v \, \delta^v, \]
for some choice of integers $k_v$. Given an element $u \in \W^j$ with $l(u) = n$, we claim that one has a pullback:
\[
\xymatrix{
\X_{ur_j}        \ar[d]^p \ar[r]^{i_{ur_j}} & \KT \ar[d]^{\pi_j} \\
\X_u^j        \ar[r]^{i_u} & \KKJ
}
\]
This follows on taking the closure of the Steinberg relation: $(\B u \B) \Gr_j = (\B u \B) \cup (\B ur_j \B)$ (see \cite{Ku} for a review of Steinberg relations). 
\noindent
Now consider the homology Serre spectral sequence for the fibration: $\KJT \rightarrow \X_{ur_j} \stackrel{p}{\rightarrow} \X_u^j$ and notice that it collapses at $E_2$. Moreover, the class $[\X_{ur_j}]$ is represented by $[\tau_j] \otimes [\X_u^j]$ at the $E_2$ stage, where we recall that $[\tau_j] \in \Hr_2(\KJT)$ is the fundamental class. Consequently, we have:
\[ 
\left<x; [\X_{ur_j}]\right>  = \left< \int \limits_{\tau_j} x; [\X_u^j] \right>, \quad x \in \Hr^*(\X_{ur_j}). \]
Taking $x = i_{ur_j}^\ast(\delta^w)$, and using the naturality of fiberwise integration, we get equalities:
\[ \left<i_{ur_j}^\ast(\delta^w); [\X_{ur_j}]\right> = \left< \int \limits_{\tau_j} i_{ur_j}^\ast(\delta^w); [\X_u^j] \right> = \sum_{v \in \W^j, \, l(v) = n} k_v \left< i_u^\ast(\delta^v); [\X_u^j] \right> = k_u. \]
But notice that the left hand side is equal to $1$ (if $ur_j = w$) and is zero otherwise. This is a restatement of the above lemma. 
\end{proof}

\medskip
\noindent
The following theorem is crucial in the study of Kac-Moody groups and their flag varieties. This structure was introduced by Kac-Peterson \cite{K1}, and studied further by Kostant-Kumar \cite{KK}. The theorem below can be shown to follow from the previous lemma, and the Bruhat exchange relations. We refer the reader to \cite{Ku} for a detailed proof. 

 \begin{thm} \label{NHR} \cite{KK} 
Define the Nil Hecke ring $\Oo$ to be the ring generated by the operators $A_i$, $i \in I$ acting on $\Hr^*(\KT)$. Then we have:
\begin{enumerate}
\item[(i)] $\Oo$ is generated by $A_i$, $i \in I$, and reltations $A_i^2 =0$, $A_i A_j A_i \cdots = 
A_j A_i A_j \cdots$ ($m_{ij}$ factors).
\item[(ii)]
$\Oo$ has an $\mbox{R}$-basis given by $\{A_w\}_{w \in W}$ where $A_w = A_{i_1} A_{i_2} \cdots A_{i_k}$ is well defined whenever $w = r_{i_1} r_{i_2} \cdots r_{i_k}$ is a reduced (or minimal length) expression.
Furthermore, 
\[
A_w (\delta^v) = \begin{cases} \delta^{vw^{-1}} \quad \mbox{if } l(w) + l(vw^{-1}) = l(v)\\
0 \quad \mbox{otherwise}
\end{cases}
\]
\end{enumerate}
\end{thm}

\begin{remark}
Notice that the algebra $\Oo$ is a deformation of the group algebra of $W$. It can be shown that it acts on $\Hr^*(\BT)$ through the operators $A_i$. The operators $A_i$ therefore generate a (twisted) algebra over $\Hr^\ast(\BT)$, which acts on $\To$-equivariant cohomology. One also has a variant of this algebra acting on the (equivariant) K-theory of the space $\KT$ \cite{Ku}, where the relation $A_i^2 = 0$ gets replaced by $A_i^2 = A_i$ (up to a unit). 
\end{remark} 

\section{Properties of $A_i$ and the action of $\W$}
\noindent
Note that $\W = \No(\To)/\To$ acts on $\To$ and consequently on $\BT$. One also has a well-defined left action of $\W$ on $\KT$ given by:
\[ w (k \To) = k w^{-1}\To, \]
where $k \To$ is a coset in $\KT$, and $w \in \W$. Recall that we have a map $\psi : \KT \rightarrow \BT$ that classified the bundle $\K \rightarrow \KT$. Up to homotopy, the map $\psi$ interpolates the action of $\W$. Consequently, we observe that the characteristic homomorphism $\psi^*: \Hr^*(\BT) \rightarrow \Hr^*(\KT)$ interpolates with the induced action of $\W$ in cohomology. In this section, we recall some fundamental results about the structure of the cohomology of Kac-Moody flag varieties as described in \cite{K1,KK,Ku}. 

\medskip
\noindent
Recall that we identified the roots $\alpha_i$, and the dual coroots $h_i^\ast$ with elements in cohomology. Using this convention, we have:

\begin{claim}\label{claim 2}
The action of $\W$ on $\Hr^2(\BT)$ and $\Hr^2(\KT)$ satisfies:
\begin{enumerate}
\item[(a)] $r_i (h_j^*) = h_j^* - \delta_{ij} \alpha_i \quad$ where $\delta_{ij}$ denotes the Kronecker
symbol and $\alpha_i \in \Hr^2(\BT)$. 
\item[(b)] $r_i (\alpha_j) = \alpha_j - a_{ij} \alpha_i$, \quad with $\alpha_i, \alpha_j \in \Hr^2(\BT)$.
\item[(c)] 
$r_i (\delta^{r_j}) = \delta^{r_j} - \delta_{ij} \alpha_i$, \quad with $\alpha_i \in \Hr^2(\KT)$. 
\end{enumerate}
\end{claim}
\begin{proof}
We may canonically identify $\Hr^2(\BT, \Ro)$ with $\gH_\Z^* \otimes \Ro$. The action of $\W$ is therefore given by $r_i(\lambda) = \lambda - \lambda(h_i) \, \alpha_i$. From this formula, part (a) and (b) follow easily. For part (c) we simply invoke Claim \ref{Euler}. \end{proof}
\noindent
Recall that the operator $A_i$ was defined via integration along the fiber for the fibration:
\[
\KIT \longrightarrow \KT \stackrel{\pi_i}{\longrightarrow} \KKI
\]
\begin{claim}\label{claim 3}
For any element $x \in \Hr^*(\KT)$, there is a unique pair of elements $y,z \in \Img \pi_i^*$ such that 
$x = \delta^{r_i} \cup y + z$. Moreover, $y = A_i(x)$ and $z = x - \delta^{r_i} \cup A_i(x)$. 
In particular, $x \in \Img \pi_i^*$ if and only if $A_i(x) =0$.
\end{claim}
\begin{proof}
The Serre spectral sequence for $\KIT \rightarrow \KT \rightarrow \KKI$ collapses at $E_2$. Since 
$\delta^{r_i}$ restricts to a generator of $\Hr^*(\KIT)$, the result follows from the multiplicative 
structure of the spectral sequence.
\end{proof}
\begin{claim}\label{claim 4}
One has the relation $\alpha_i \cup A_i (x) = x - r_i (x)$ for $x \in \Hr^*(\KT)$.
\end{claim}
\begin{proof}
Write $x$ as $x = \delta^{r_i} \cup A_i(x) + z$ with $z \in \Img \pi_i^*$. Since $r_i$ acts on $\KT$ through right multiplication by an element in $\K_i$, it fixes the image
of $\pi_i^*$. Thus $r_i (x) = (\delta^{r_i} - \alpha_i ) \cup A_i(x) + z $ using \ref{claim 2}. 
The difference of the above two equations gives the required result.
\end{proof}
\begin{claim} \label{claim 5}
If $2 \in \Ro$ is not a zero divisor, then $\alpha_i$ is not a zero divisor on the image of $\pi_i^*$. In particular, under this assumption, $A_i(x)$ is the unique
element in the image of $\pi_i^*$ satisfying the equality given in \ref{claim 4}.

\end{claim}
\begin{proof}
Assume $\alpha_i \cup y =0$ for some $y \in \Img \pi_i^*$. Note that $\alpha_i = 2 \delta^{r_i} + z$ for some 
$z \in \Img \pi_i^*$ using \ref{Euler}. Thus we have $ \delta^{r_i} \cup 2y + z \cup y =0$. From \ref{claim 3} we get 
$2y=0$ which implies $y =0$ by assumption.
\end{proof}

\begin{remark} \label{Steenrod} The same proof as above shows that 
\ref{claim 3}, \ref{claim 4} and \ref{claim 5} hold for $\Hr^*(\BT)$  once we replace $\delta^{r_i}$ by $h_i^*$, and $\pi_i^*$ by $\theta_i^*$.
\end{remark}
\begin{thm} \label{A of a product}
$A_i(u \cup v) = A_i(u) \cup r_i(v) + u \cup A_i(v)$, where $u,v$ are arbitrary homogeneous elements in $\Hr^*(\KT)$ or $\Hr^*(\BT)$. 
\end{thm}
\begin{proof}
Since the operators $A_i$ are natural with respect to ring homomorphisms, it is sufficient to prove the theorem for $\Ro = \Z$. Note that:
\[
\begin{split}
\alpha_i \cup (A_i(u) \cup r_i(v) + u \cup A_i(v)) &= (u - r_i (u)) \cup r_i(v) + u \cup (v - r_i (v))\\
& =u \cup v - r_i (u \cup v)\\
& = \alpha_i \cup A_i (u \cup v)
\end{split}
\]
Thus, by \ref{claim 5}, we will be done if we can show that $ x = A_i(u) \cup r_i(v) + u \cup A_i(v)
\in \Img \pi_i^*$. It suffices to show that $\alpha_i \cup A_i(x) =0$, but we have $\alpha_i \cup A_i (x) = x - r_i (x) = -\alpha_i \cup A_i(u) \cup A_i(v) + \alpha_i \cup A_i(u) \cup A_i(v) =0$ and so we are done.
\end{proof}
\begin{thm}
Let $\Oo_J \subseteq \Oo$ be the subalgebra generated by the operators $A_j, j \in J \subseteq I$. Then 
$\Hr^*(\K/\K_J)= \Hr^*(\KT)^{\Oo_J}$, where $\Hr^*(\KT)^{\Oo_J}$ denotes all elements annihilated by $\Oo_J$ and $\Hr^*(\K/\K_J)$ is identified with its image in $\Hr^*(\KT)$ via $\pi_J^*$.
\end{thm}
\begin{proof}
The image of $\pi_J^*$ is a free $\Ro$-module on the Schubert basis $\delta^w$, $w \in \W^J$. But $w \in \W^J$ if
and only if $l(wr_j) > l(w)$, for all $j \in J$. The result follows from the formula for the action of $A_j$ on the Schubert basis \ref{Schubert}.
\end{proof}
\begin{corr}
If $2 \in \mbox{R}$ is not a zero divisor, then $\Hr^*(\K/\K_J) = \Hr^*(\KT)^{\W_J}$ is the submodule of $\W_J$-invariant elements for $J \subseteq I$. In particular $\Hr^*(\KT)^{\W} = \Hr^0(\KT) = \Ro$.
\end{corr}
\begin{proof}
$\Hr^*(\K/\K_J) = \Hr^*(\KT)^{\Oo_J} = \Hr^*(\KT)^{\W_J}$ using \ref{claim 4} and \ref{claim 5}. 
\end{proof}

\medskip
\begin{remark}
Let $\Ro = \F_2$ and consider the example of the Kac-Moody group corresponding to the Cartan matrix $\left( \begin{matrix} 2 & -2 \\ -2 & 2 \end{matrix} \right)$. Then both the simple roots are zero, and so \ref{claim 4} shows that $\Hr^*(\KT; \F_2)^{\W} = \Hr^*(\KT; \F_2)$. This example shows that the assumption that $2$ is not a zero-divisor is necessary. 
\end{remark}

\section{The cohomology of $\K$}
\noindent
In this section we recall some results from \cite{K}. Recall the characteristic homomorphism $\psi^*: \Hr^*(\BT) \rightarrow \Hr^*(\KT)$. We begin by identifying the kernel of $\psi^*$. 
Let $\Hr^*(\BT)^+$ denote elements of $\Hr^*(\BT)$ in the kernel of the augmentation to $\Hr^0(\BT)$. Similarly, denote $\Hr^\ast(\KT)^+$ to be the elements in the kernel of the augmentation to $\Hr^0(\KT)$. 

\medskip
\noindent
Let $J \subset \Hr^*(\BT)$ be defined as: $J = \{u \in \Hr^*(\BT) \, | \, A_{i_1} A_{i_2} \cdots A_{i_k} (u) \in \Hr^*(\BT)^+ \, \forall i_1,\ldots,i_k \}$. We have:
\begin{thm} \cite{K}
$J = \Ker \psi^*$.
\end{thm}
\begin{proof}
Notice first using \ref{NHR} that $z \in \Hr^*(\KT)^+$ is nonzero if and only if there exists a sequence $i_1 \ldots i_k$ such
that $0 \neq A_{i_1} \cdots A_{i_k} (z) \in \Hr^0(\KT)$. Now since $\psi^*$ is an isomorphism in degree $0$, the statement $A_{i_1} \cdots A_{i_k} (u) \in \Hr^*(\BT)^+ \, \forall i_1 \ldots i_k$ is equivalent to the statement $\psi^*A_{i_1} \cdots A_{i_k} (u) \in \Hr^*(\KT)^+ \, \forall i_1 \ldots i_k$, i.e. $A_{i_1} \cdots A_{i_k}
\psi^*(u) \in \Hr^*(\KT)^+ \, \forall i_1 \ldots i_k$ which is equivalent to $\psi^*(u) =0$.
\end{proof}
\begin{remark} 
The ideal $J$ is known as the ideal of {\em Generalized Invariants} of $\W$. It has been studied in detail and appears to be of independent interest \cite{KP4}. 
\end{remark}

\medskip
\noindent 
For the rest of this section we work with field coefficients, we call our field $\F$.  

\begin{lemma} \cite{K}
$J/J^2$ is a free $\Hr^*(\BT)$-module.
\end{lemma}
\begin{proof}
Let $y_1,y_2,\ldots$ be a set of homogeneous elements of $J$ such that $y_1, y_2, \ldots$ form a minimal set of generators of $J/J^2$ as an $\Hr^*(\BT)$-module ordered by increasing degree: $0 < \operatorname{deg} (y_1) \leq \operatorname{deg}(y_2) \leq \operatorname{deg} (y_3) \ldots$. Now recall that the characteristic homomorphism $\psi^\ast$ intertwines the operators $A_i$ as well as the action of $\W$. Therefore, \ref{A of a product} implies that the operators $A_i$ preserve $J^2$ and that:
\[
A_i (y_k) = \sum_{j<k} r_j y_j + J^2, \qquad (*)
\]
with $r_j \in \Hr^*(\BT)$. Let $\sum_{j \leq k} s_j  y_j \in 
J^2$ be some homogeneous relation. We can assume $s_k \notin J$. Choose a sequence $i_1 \ldots i_k$ such that $0 \neq A_{i_1} \cdots A_{i_k} (s_k) \in \F$.
Applying $A_{i_1} \cdots A_{i_k} $ to the relation and using $(*)$ repeatedly we notice that:
\[
y_k = \sum \limits_{j<k} t_j y_j + J^2 \quad t_j \in \Hr^*(\BT)
\]
which is a contradiction to the minimality of the set of generators.
\end{proof}
\begin{thm} \cite{K}
$J$ is generated by a regular sequence $J = \langle \sigma_1, \sigma_2, \ldots, \sigma_r \rangle$, where $r \leq$ rank$(\To)$. In particular, $\Hr^\ast(\BT)$ is a free module over the polynomial algebra generated by the classes $\sigma_i$. 
\end{thm}
\begin{proof}
First notice that any regualr sequence in $\Hr^*(\BT)$ must have length $\leq$ rank$(\To)$. So it remains to show that $J$ is generated by a regular sequence. This follows from a theorem of Vasconcelos \cite{V} which says that for a graded algebra $\mbox{A}$ of finite global dimension a homogeneous ideal $J \subseteq \mbox{A}$ is generated by a regular sequence if and only if $J/J^2$ is a free $\mbox{A}/J$-module.
\end{proof}
\begin{remark}
If $\F$ has positive characteristic $p$, the length of this sequence is exactly rank$(\To)$. This can be seen as follows: First notice that if $\lambda \in \Hr^*(\BT)^{\W}$, then its $p$-th power $\lambda^p$ is annhilated by all the operators $A_i$. Now notice that $\Hr^*(\BT)^{\W}$ contains the Dickson invariants \cite{D} defined as the invariants with respect to the action of $\mbox{GL}_n(\F)$, where $n$ is the rank of $\To$. The Dickson invariants (or their $p$-th powers) form a regular sequence of length $n$. Hence $J$ contains a regular sequence of maximal length and must therefore itself be generated by a sequence of maximal length. 
\end{remark}

\medskip
\noindent
Define a subring $S$ of $\Hr^*(\KT)$ via $S= \Img \{\psi^*: \Hr^*(\BT) \rightarrow \Hr^*(\KT)\} \iso \Hr^*(\BT)/J$.
Notice that:
\begin{enumerate}
\item[a)] $S$ is the subring  of $\Hr^*(\KT)$ generated by $\Hr^2(\KT)$ since the map $\KT \rightarrow \BT$ is $2$-connected.
\item[b)] If $\F$ has characteristic $p>0$ then $S$ is a finite dimensional vector space over $\F$. In fact, since it is a complete intersection ring, it has the structure of of a Poincar\'{e} duality algebra, though we will have no use for this fact in this document. 
\item[c)] Let $2d_i$ be the degree of $\sigma_i$, and let $n$ be the rank of $\To$, then the Poincar\'e series of $S$ is: 
\[ P_t(S) = \frac{\prod_{i=1}^r (1-t^{2d_i})}{(1-t^2)^n}, \] 
as can be easily seen using the fact that $\Hr^\ast(\BT)$ is a free module over the polynomial algebra generated by the classes $\sigma_i$. 

\end{enumerate}
\begin{remark}
The facts a) and c) are true in arbitrary characteristic, but fact b) requires the characteristic to be positive. Indeed, examples in rank two (see section 10), show that the regular sequence generating $J$ has strictly less length than the rank, and so $S$ cannot be finite dimensional. 
\end{remark}
\begin{thm} \cite{K}
$\Hr^*(\KT)$ is a free $S$-module.
\end{thm}
\begin{proof}
Proceed as before. Let $y_1,y_2 \ldots$ be a set of homogeneous elements of $\Hr^*(\KT)$ so that $y_1,y_2 \ldots$ form an $\F$-basis of $\Hr^*(\KT) \otimes_S \F$, and $0 = \operatorname{deg} (y_1) \leq \operatorname{deg}(y_2) \leq \operatorname{deg} (y_3) \ldots$. It is clear that the $y_i$ generate $\Hr^*(\KT)$ as an $S$-module. Note that the operators $A_i$ preserve $S$ for all $i$ and 
\[
A_i(y_k) = \sum_{j<k} r_j y_j; \quad r_j \in S \qquad (*)
\]
Let $\sum_{j\leq k} s_j y_j =0$ be some homogeneous relation in $\Hr^*(\KT)$ and assume that $s_k \neq 0$.
Choose a sequence $i_1 \ldots i_k$ sucht that $0 \neq A_{i_1} \cdots A_{i_k} (s_k) \in \F$. Applying 
$A_{i_1} \cdots A_{i_k}$ to the relation and using $(*)$ repeatedly we get a contradiction.
\end{proof}

\medskip
\noindent 
We proceed now to study some standard fibrations related to $\K$. Recall that given a fibration $F \rightarrow E \rightarrow B$, 
with $B$ simply connected, the Eilenberg-Moore spectral sequence is a second quadrant cohomological spectral sequence of graded algebras \cite{JM}. The $E_2$ term is given by $E_2^{p,q} = \Tor^{-p,q}_{\Hr^*(B)}(\F, \Hr^*(E))$, and it converges to $\Hr^*(F,\F)$ for a field $\F$. As an easy consequence of the results we have so far, we recover the result of Kac \cite{K}: 
\begin{thm} 
The Eilenberg-Moore spectral sequence for $ \K \stackrel{\pi}{\rightarrow} \KT \stackrel{\psi}{\rightarrow}
\BT$ collapses at $E_2$. Furthermore, $\Hr^*(\K,\F)$ is free over the image of $\pi^*$, and one has a short exact sequence of algebras:
\[ 1 \rightarrow  \Hr^*(\KT,\F)\otimes_S \F \stackrel{\pi^*} \longrightarrow \Hr^*(\K,\F) \longrightarrow \Lambda(x_1,\ldots x_r) \rightarrow 1. \]
\end{thm}
\begin{proof} 
We apply the spectral sequence to the fibration $\K \stackrel{\pi} \rightarrow \KT \stackrel{\psi} \rightarrow \BT$ in order to compute $\Hr^*(\K)$. The $E_2$-term is given by:
\[ E_2^{p,q} = \Tor^{-p,q}_{\Hr^*(\BT)}(\F, \Hr^*(\KT)). \]
From the previous theorem, we may write $\Hr^*(\KT)$ as the $S$-module $S \otimes_\F \Hr^*(\KT) \otimes_S \F$. It follows that 
\[ E_2^{p,q} = \Tor^{-p,q}_{\Hr^*(\BT)}(\F,S) \otimes_\F \Hr^*(\KT) \otimes_S \F = \Lambda(x_1,x_2,\ldots, x_r) \otimes_\F \Hr^*(\KT) \otimes_S \F, \]
where $\Lambda(x_1,x_2,\ldots,x_r)$ denotes an exterior algebra on classes $x_i$ of homogeneous bidegree ($-1, |\sigma_i|$). Due to degree reasons, this spectral sequence collapes. Consider the subring $\Hr^*(\KT) \otimes_S \F$. Since it is in bidegree ($0,\ast$), it can be identified via the edge homomorphism with the image of the map $\pi^\ast : \Hr^*(\KT) \rightarrow \Hr^\ast(\K)$. This proves the above theorem. 
\end{proof}

\begin{remark}
It is natural to ask if $\Hr^*(\KT)$ and $\Hr^*(\K)$ are finitely generated algebras. If the characteristic of the field $\F$ nonzero, then we will use results from the next section to show that both these algebras are not finitely generated, unless $\K$ is a compact Lie group. In characteristic zero, this question remains open. 
\end{remark}

\section{$\Hr^*(\KT)$ as a module over the Steenrod algebra $\mathcal{A}_p$}
\noindent
In this section, we study the structure of the cohomology of $\K$ and $\KT$ as modules over the mod $p$ Steenrod algebra $\mathcal{A}_p$. 
\begin{thm} \cite{K1}
Let $p$ a prime and let $\F =\F_p$. Let $\mathcal{P} = \sum P^i$ be the total Steenrod operation ($P^i = 
\operatorname{Sq}^i$ if $p=2$). Then $A_i (\mathcal{P} (x)) = (1 + \alpha_i^{p-1}) \mathcal{P} (A_i(x))$ for any element $x \in \Hr^*(\KT)$.
\end{thm}
\begin{proof}
First recall by remark \ref{Steenrod} that one has a relation in $\Hr^\ast(\BT,\Z)$:
\[ \alpha_i \cup A_i(x) = x - r_1(x). \]
Now $\alpha_i \in \Hr^2(\BT,\Z)$ is not a zero divisor, and so we may use the $\W$ action on $h_i^\ast$ given in \ref{claim 2} to obtain the well defined formula over $\Z$: 
\[ A_i( (h_i^*)^p ) = \frac{(h_i^*)^p - (h_i^* -\alpha_i)^p}{\alpha_i} \equiv (-\alpha_i)^{p-1} \equiv (\alpha_i)^{p-1} \operatorname{mod} p. \]
Applying the characteristic homomorphism $\psi^*$ and reducing  modulo $p$ we get the formula:
$A_i((\delta^{r_i})^p ) = \alpha_i^{p-1}$. Let $x \in \Hr^*(\KT)$. We may express it as $x = \delta^{r_i} \cup A_i(x)  + z$. Then
$\mathcal{P} (x) = (\delta^{r_i} + (\delta^{r_i})^p) \cup \mathcal{P} (A_i (x)) + \mathcal{P} (z)$. And thus
$A_i (\mathcal{P}(x)) = (1 + \alpha_i^{p-1}) \cup \mathcal{P} (A_i (x))$ using the fact that $A_i$ is a map of $\Hr^*(\KKI)$-modules.
\end{proof}

\noindent
Now let $J = \langle \sigma_1, \sigma_2, \ldots, \sigma_r \rangle$ be the ideal of generalized invariants. Let $2d_i$ denote the degree of the element $\sigma_i$. The degree of the top class in $S$ is given by $2m = 2 \sum (d_i -1)$. We will use the above theorem to show that $\Hr^*(\KT)$ (and therefore $\Hr^*(\K)$) is locally finite as a module over $\mathcal{A}_p$ (i.e. $\mathcal{A}_p (z)$ is a finite dimensional vector space for all $z \in \Hr^*(\KT)$).

\medskip
\noindent 
For any homogeneous subset $X\subseteq \Hr^*(\KT)$ define $d(X)\leq \infty$ to be the highest degree
of any homogeneous element in $X$. For $z \in \Hr^*(\KT)$. let $M(z)$ be the $S$-module given by the span of elements of the form $s \cup a(z)$, $s \in S$, $a \in \mathcal{A}_p$. Note that $M(z)$ is a n 
$\mathcal{A}_p$-submodule of $\Hr^*(\KT)$. Let $d(z)$ denote $d(M(z))$.
\begin{thm}
If $z$ is an element of positive homogeneous degree $2k$, then $d(z) \leq 2 k (m+1) -2$.
\end{thm}
\begin{proof}
We work by induction on the degree of $z$. Since $S$ is the subring of $\Hr^*(\KT)$ generated by elements of degree 2, we are done for $k=1$. Now let $z$ be any element of homogeneous degree $2 k + 2$. Let $x = \sum_\mu s_\mu \cup P^\mu (z)$ be a homogeneous element of $M(z)$ where $\mu$ ranges over finite sequences
of positive integers $i_1 \ldots i_s$ and $P^\mu = P^{i_1} \cdots P^{i_s}$. By repeated application of the previous theorem we notice that $A_i (P^\mu (z)) \in M(A_i(z))$ for any $i \in I$. By induction, $P^\mu(z)$ can have
degree at most $2 k (m+1) -2 +2 = 2k(m+1)$. Thus $x$ has degree at most $2 k (m+1) + 2m$. 
Hence $d(z) = d(M(z))\leq 2 m + 2 k(m+1) = 2(k+1)(m+1) -2$ and we are done.
\end{proof}
\begin{corr}
$\Hr^*(\KT)$ and $\Hr^*(\K)$ are locally finite as modules over $\mathcal{A}_p$. In particular, over a field $\F$ of positive characteristic $\Hr^*(\KT)$ is finitely generated if and only if it is finite dimensional (i.e. if and only if $\K$ is a compact Lie group). 
\end{corr}
\begin{proof} The first part of the statement follows from the previous theorem . For the second part, notice that for $\Hr^*(\KT)$ to be infinite dimensional and finitely generated, there must exist an element $\lambda$ with arbitrary large nonzero powers. This is impossible since $\mathcal{A}_p$ acts locally finitely on $\lambda$. The same argument works for $\Hr^*(\K)$. 
\end{proof}
\section{The Hopf algebra structure of $\Hr^*(\K)$}
\noindent
Recall the extension of algebras:
\[
1 \rightarrow \Hr^*(\KT)\otimes_S \F \stackrel{\pi^*} \longrightarrow \Hr^*(\K) \longrightarrow \Lambda(x_1,\ldots,x_r) \rightarrow 1.
\]
We will show that this is actually an extension of Hopf algebras. It will be sufficient to show that 
$\Hr^*(\KT) \otimes_S \F$ is a sub coalgebra of $\Hr^*(\K)$. 

\noindent We define a coalgebra structure on $\Hr^*(\KT)$ (introduced by D. Peterson) via:
\[
\Delta (\delta^w) = \sum \limits_{uv = w} \delta^u \tens \delta^v
\]
where the sum runs over all reduced expressions of $w$. Recall that $\Hr^*(\KT) \otimes_S \F$ maps isomorphically to the image of $\pi^*$ in $\Hr^*(\K)$. Our main theorem of this section will state that $\pi^*: \Hr^*(\KT) \rightarrow \Hr^*(\K)$ is a map of coalgebras. The first step towards this goal is the construction of equivariant Schubert classes.

\medskip
\noindent 
Let $\Hr^*_{\To}(\KT) = \Hr^*(\ET \times_{\To} (\KT))$ denote the equivariant cohomology 
of $\KT$. For the moment, we allow coeffcients in any ring. Define homomorphisms $E_w$ for 
$w \in \W$ by: 
\[
E_w: \Hr^*_{\To}(\KT) \stackrel{i_w^*}{\longrightarrow} \Hr^*_{\To}(\X_w) \stackrel{\int_{[\X_w]}}{\longrightarrow}
\Hr^*_{\To}(pt)
\]
where $\int_{[\X_w]}$ denotes integration over the fiber for the oriented fibration:
\[
\X_w \longrightarrow \ET \times_{\To} \X_w \longrightarrow \BT
\]
Note that these are homomorphisms of $\Hr^*_{\To}(pt)$-modules.

\begin{claim}
There exists a unique basis $\{ \delta^w_{\To} \}_{w \in \W}$ of $\Hr^*_{\To}(\KT)$ over $\Hr^*_{\To}(pt)$ with the property:
\[
E_v(\delta^w_{\To}) = \begin{cases}
1 \quad \text{if} \quad v =w\\
0 \quad \text{otherwise}
\end{cases}.
\]
\end{claim} 
\begin{proof}
Uniqueness will follow easily once we have existance. We proceed by induction on $l(w)$. For 
$w =1$, let $\delta^w_{\To} =1$. Assume we are done defining $\delta^w_{\To}$ for $l(w)<k$.
Using the naturality of the pullback:
\[
\xymatrix{
\X_w       \ar[d] \ar[r] & \ET \times_{\To} \X_w \ar[d] \\
pt        \ar[r] & \BT
}
\]
we see that $\delta^w_{\To}$ restricts to $\delta^w$ under the restriction map $\iota^* : \Hr^*_{\To}(\KT) \rightarrow \Hr^*(\KT)$.

\noindent
For $w \in \W$ such that $l(w) =k$, let $x_w \in \Hr^{2k}_{\To}(\KT)$ be any element that restricts
to $\delta^w$ under the (surjecive) map $\iota^*$. For degree reasons, we see that:
\begin{alignat*}{2}
E_w(x_w) &= 1 &\quad & \text{and}\\
E_v(x_w) & =0 & \quad & \text{if} \quad l(v) \geq k, w \neq v
\end{alignat*}
Now $\delta^w_{\To}$ can be defined as $\delta^w_{\To} = x_w - \sum \limits_{l(v) <k} E_v(x_w) 
\delta_{\To}^v$, and we are done by induction.
\end{proof}
\noindent
One has a similar fibrations:
\[ \KT \longrightarrow \K \times_{\To} (\KT) \stackrel{p} \longrightarrow \KT, \qquad X_w \longrightarrow \K \times_{\To} X_w \stackrel{p} \longrightarrow \KT, \]
and we can define homomorphisms:
\[
F_w: \Hr^*(\K \times_{\To} (\KT)) \stackrel{i_w^*}{\rightarrow} \Hr^*(\K \times_{\To} X_w)
\stackrel{\int_{[\X_w]}}{\rightarrow} \Hr^*(\KT).
\]
Note that $F_w$ is a homomorphism of $\Hr^*(\KT)$-modules where $\Hr^*(\K \times_{\To}(\KT))$ is 
viewed as a $\Hr^*(\KT)$-module via $p^*$. 

\medskip
\noindent
Consider the pullback obtained by performing the associated bundle over $\psi$ for the left $\To$-action on $\KT$:
\[
\xymatrix{
\K \times_{\To}(\KT)   \ar[d]^p \ar[r]^{\Psi} & \ET \times_{\To} (\KT) \ar[d] \\
\KT        \ar[r]^{\psi} & \BT
}
\]
\begin{claim}
Let $\sigma^w = \Psi^* (\delta^w_{\To})$. Then the set $\{ \sigma^w\}_{w \in W}$ is the unique
basis of $\Hr^*(\K \times_{\To} (\KT))$ as an $\Hr^*(\KT)$-module with the property:
\[
F_v(\sigma^w) = \begin{cases}
1 \quad &\text{if} \quad v = w\\
0 \quad &\text{otherwise}
\end{cases}
\]
\end{claim}
\begin{proof}
\noindent
The proof follows from the previous claim and naturality with respect to pullback for the above diagram.
\end{proof}
\noindent
We will now need some technical results which will be the content of the next few lemmas.
First consider the bundle over $\C P^1$ given by $\K_i \rightarrow \KIT$. Fix sections $\Y_i$ over the open cell $\C \subset \C P^1$ for each $i \in I$. Consider the following subspace of $\K$:
\[
\Zo_u = \bigcup \limits_{u' \leq u} \Y_{u'} \cdot \To.
\]
where $\Y_u$ is the space $\Y_u = \Y_{i_1} \times \ldots \times \Y_{i_k}$ with $u = r_{i_1} \ldots r_{i_k}$ is a reduced expression. For $u,v \in \W$, we can form the space $\X_{u,v} = \Zo_u \times_{\To} \X_v $. The space $\X_{u,v}$ has a decomposition:
\[
\X_{u,v} = \bigcup \limits_{u' \leq u, v'\leq v} (\Y_{u'} \cdot \To) \times_{\To} \Y_{v'} = 
\bigcup \limits_{u' \leq u , v' \leq v }\Y_{u'} \times \Y_{v'}.
\]
The CW decomposition of $\KT$ implies that the space $\K \times_{\To} \KT$ admits a CW decomposition with even open cells given by the spaces $\Y_u \times \Y_v$. Let $\delta^{u,v} \in \Hr^{2l(u) + 2 l(v)} (\X_{u,v})$ be the class linear dual to the cell: $[\Y_u \times \Y_v]$.

\medskip
\noindent
Let 
$\mu: \K \times_{\To} (\KT) \rightarrow \KT$ denote the left action of $\K$ on $\KT$. One has
 an induced map:
\[
\mu_{u,v}: \X_{u,v} = \Zo_u \times_{\To} \X_v \longrightarrow \K \times_{\To} (\KT)
\stackrel{\mu}{\longrightarrow} \KT.
\]

\begin{lemma} \label{lemma 1}
If $w \in \W$ such that $l(w) = l(u) + l(v)$ then:
\[
\mu_{u,v}^* (\delta^w) = 
\begin{cases} 
\delta^{u,v} \quad & \text{if} \quad uv =w\\
0 \quad &\text{otherwise}
\end{cases}
\]
\end{lemma}
\begin{proof}
Recall from the CW-decomposition of $\KT$ that we have a unique factorization of cells in $\KT$ \cite{Ku} given by: 
$\Y_u \cdot \Y_v = \Y_{uv}$ if $uv$ is reduced. Otherwise, $\Y_u \cdot \Y_v$ factors through cells
of lower dimension. The result follows.
\end{proof}
\noindent
Now consider the oriented fibration
\[
\X_v \longrightarrow \Zo_u \times_{\To} \X_v \stackrel{p_{u,v}}{\longrightarrow} \X_u
\]
\begin{lemma}\label{lemma 2}
\[
\int_{[\X_v]} \delta^{u,v} = \delta^u
\]
\end{lemma}
\begin{proof}
The proof is obvious since the Serre sepctral sequence for the above filtration 
collapses at $E_2$ and $\delta^{u,v}$ is represented by $\delta^u \tens \delta^v$.
\end{proof}

\noindent
Recall that $\Hr^\ast(\K \times_{\To} \KT)$ is a free $\Hr^\ast(\KT)$-module via $p^*$, on the basis $\sigma^v$. We have:

\begin{lemma}\label{lemma 3}
Under the action map $\mu: \K \times_{\To} (\KT) \longrightarrow \KT$, we have the equality:
\[
\mu^*(\delta^w) = \sum \limits_{uv=w} \delta^u \cup \sigma^v
\]
where the sums runs over all reduced expressions of $w$. 
\end{lemma}
\begin{proof}
Let $l(w) =k$. We can write 
\[ \mu^*(\delta^w) = \sum \limits_{l(u) + l(v) =k} a_{u,v} \, \delta^u \cup \sigma^v\] 
where $a_{u,v}$ are elements in the coefficient ring. We use the operators $F_v$ to isolate
them:
\[
F_v \, \mu^*(\delta^w) = \sum \limits_{l(u) = k - l(v)} a_{u,v} \, \delta^u
\]
and therefore $i_u^*\, F_v \, \mu^*(\delta^w) = a_{u,v} \, \delta^u$ where $i_u^\ast: \Hr^\ast(\KT) \rightarrow \Hr^\ast(\X_u)$ is induced by the
inclusion. Now consider the commutative diagram:
\[
\xymatrix{
\Zo_u \times_{\To} \X_v   \ar[d]^{\mu_{u,v}} \ar[r] & \K \times_{\To} \X_v \ar[r]^{i_v} \ar[d]^{\mu} & \K \times_{\To} (\KT)  \ar[d]^{\mu} \\
\X_{uv}        \ar[r] & \KT \ar[r]^{=} & \KT
}
\]
We have:
\[
i_u^* \, F_v \, \mu^*(\delta^w) = \int_{[X_v]} \mu^*_{u,v}(\delta^w) = \begin{cases}
\delta^u \quad & \text{if} \quad uv=w\\
0 \quad &\text{otherwise}
\end{cases}
\]
using \ref{lemma 1} and \ref{lemma 2}, so we are done.
\end{proof}
\noindent
Let $m : \K \times \K \rightarrow \K$ denote the multiplication map. Consider the commutative diagram:
\[
\xymatrix{
\K \times \K   \ar[d]^m \ar[r]^{\tilde{\pi}\quad } & \K \times_{\To} (\KT) \ar[d]^{\mu} \\
\K        \ar[r]^{\pi} & \KT
}
\]
Let $\pi_1$, $\pi_2$ be defined
as the projections onto the first and second factor respectively:
\begin{lemma} \label{lemma 4}
In cohomology, $\tilde{\pi}$ is given by
\begin{enumerate}
\item[a)] $\tilde{\pi}^* (\sigma^v) = \pi_2^* \circ \pi^* (\delta^v)$
\item[b)] $ \tilde{\pi}^*(\delta^u) = \pi_1^* \circ \pi^*(\delta^u)$
\end{enumerate}
\end{lemma}
\begin{proof}
For a), recall that by definition $\sigma^v = \Psi^*(\delta_{\To}^v)$. Now notice that
we have a commutative diagram:
\[
\xymatrix{
\K \times \K   \ar[d]^{\tilde{\pi}} \ar[r]^{f} & \ET \times (\KT) \ar[d]^{g} \\
\K \times_{\To} (\KT)       \ar[r]^{\Psi} & \ET \times_{\To} (\KT)
}
\]
Note that $\ET \times (\KT) \stackrel{\sim}{\rightarrow} \KT$ since $\ET$ is contractible. Under this
identification $f = \pi \circ \pi_2$ and $g$ is the inclusion of $\KT$ in $\ET \times_{\To} (\KT)$.
Since $g^*(\delta^v_{\To}) = \delta^v$, we get:
\[ \tilde{\pi}^* (\sigma^v) = \tilde{\pi}^* \circ \Psi^* (\delta_{\To}^v) = \pi_2^* \circ \pi^* (\delta^v). \]

\noindent Part b) follows from the commutative diagram:
\[
\xymatrix{
\K \times \K   \ar[d]^{\pi_1} \ar[r]^{\tilde{\pi}} & \K \times_{\To} (\KT) \ar[d]^{p} \\
\K       \ar[r]^{\pi} & \KT.
}
\]
\end{proof}
\noindent
We now prove the main theorem of this section:
\begin{thm}
The map $\pi^*: \Hr^*(\KT) \rightarrow \Hr^*(\K)$ is a map of coalgebras, where cohomology is 
taken with coefficient in a field.
\end{thm}
\begin{proof}
One recalls the commutative diagram:
\[
\xymatrix{
\K \times \K   \ar[d]^m \ar[r]^{\tilde{\pi}} & \K \times_{\To} (\KT) \ar[d]^{\mu} \\
\K        \ar[r]^{\pi} & \KT
}
\]
Now one invokes \ref{lemma 3} and \ref{lemma 4}, to get the required equality:
\[ m^*\pi^*(\delta^w) = 
\sum \limits_{uv=w} \pi_1^*(\pi^* (\delta^u)) \cup \pi_2^*(\pi^* (\delta^v)), \]
where the sum is begin taken over all reduced expressions. This is exactly the statement of the theorem.
\end{proof}
\begin{remark} 
Note that \ref{lemma 1}, \ref{lemma 2}, \ref{lemma 3} and \ref{lemma 4} made no assumption
on the coefficient ring. For an arbitrary ring $\Ro$, define $A_{\Ro} = \Img (\pi^* \Hr^*(\KT; \Ro) \rightarrow \Hr^*(\K;\Ro))$. Note that for  a field $\F$, we have shown that $A_\F = \Hr^*(\KT; \F) \otimes_S \F$. It follows easily that \cite{K1}:
\begin{enumerate}
\item[1)] $A_\F = A_\Z \otimes_\Z \F$,  
\item[2)] $A_\Z$ has the structure of a Hopf algebra. 

\end{enumerate}
Thus $A_\Z$ is an integral lift of the Hopf algebras $A_\F$.
\end{remark}

\section{On the ring structure of $\Hr_*(\K)$}

\noindent
We fix a coefficient field $\F$ of positive characteristic throughout this section. It is our aim to show that the Pontrjagin ring: $\Hr_*(\K) = \Hr_*(\K,\F)$ is a finitely generated algebra. We begin with some preliminary lemmas:

\begin{lemma} The left action of $\Hr_*(\K)$ on $\Hr_*(\KT)$ factors through the projection $\Hr_*(\K) \rightarrow A_{\F}^*$, where $A_{\F}^\ast$ is the dual Hopf algebra of $A_{\F}$. Moreover, $\Hr_*(\KT)$ is a finitely generated free left $A_{\F}^*$- module. 
\end{lemma}
\begin{proof}
The first part of the above lemma follows by dualizing the result of \ref{lemma 3}. The theory of Hopf-algebras can now be used to show that $\Hr_*(\KT)$ is a free $A_{\F}^*$-module \cite{MM}(Thm. 4.4). Now since $\Hr^\ast(\KT)$ is a free module over the finite dimensional algebra $S$, with a basis given by $A_{\F}$, we can compute Poincar\'e series to see that $\Hr_\ast(\KT)$ must have a finite basis over $A_{\F}^\ast$. 

\smallskip
\noindent
For the benefit of the reader, we also provide an alternate proof: Recall that for a topological group $G$, and a principal $G$-bundle $G \rightarrow E \rightarrow B$,
one has a natural homological (Bar) spectral sequence of coalgebras converging to $\Hr_*(B)$ with $E_2$-term given by $E_2^{p,q} = \Tor_{\Hr_*(G)}(\F, \Hr_*(E))$ \cite{JM}. 
Consider a pair of pullbacks of principal $\To$-bundles:
\[
\xymatrix{
\K \times \K   \ar[d]^{Id \times \pi} \ar[r]^{m} & \K  \ar[d]^{\pi} \ar[r] & \ET \ar[d] \\
\K \times (\KT)       \ar[r]^{\mu} & \KT \ar[r]^{\psi} & \BT
}
\]
This induces an action of $\Hr_*(\K)$ on the Bar spectral sequence of converging to $\Hr_*(\KT)$:
\[ \Hr_*(\K) \otimes \Tor_{\Hr_*(\To)}(\F, \Hr_*(\K)) \longrightarrow \Tor_{\Hr_*(\To)}(\F, \Hr_*(\K)). \]
Since $\K$ is simply connected, the inclusion of the maximal torus $\To \subset \K$ is null homotopic. Consequently, we have:
\[ \Tor_{\Hr_*(\To)}(\F, \Hr_*(\K)) = \Hr_*(\K) \otimes_\F \Tor_{\Hr_*(\To)}(\F,\F). \]
Differentials in this spectral sequence must must annhilate the piece corresponding to the dual exterior algebra: $\Lambda^* (x_1, \ldots, x_r) \subseteq \Hr_*(\K)$. This dual algebra is itself an exterior algebra. The generators of this exterior algebra must therefore be targets of differentials originating on elements indecomposable under the $\Hr_*(\K)$-action. Now we may write: 
\[ \Tor_{\Hr_*(\To)}(\F, \F) = \Gamma(y_1,\ldots,y_r) = S^* \otimes_\F \Gamma(\tau_1,\ldots,\tau_r) \]
where $\Gamma(y_1,\ldots,y_r)$ denotes the dual of a polynomial algebra. This coalgebra is bigraded by giving $y_i$ bidegree $(1,1)$. The element $\tau_i$ is an element of bidegree $(1,|x_i|)$, the vectorspace $S^*$ is dual to $S$, and is detected in $\Hr_*(\BT)$. It follows that the generators $\tau_i$ must hit a set of generators of the dual exterior algebra $\Lambda^*(x_1,\ldots, x_r)$ in the spectral sequence. Consequently, the $E_{\infty}$ term of the spectral sequence is a free left module over $A_\F^* = \Hr_*(\K) \otimes_{\Lambda^*(x_1,\ldots, x_r)} \F$, with a basis given by the finite dimensional vector space $S^*$. The result follows from an easy filtration argument.  
\end{proof} 

\noindent
We now need the following general lemma:

\begin{lemma} Let $A$ be a (not necessarily commutative) finitely generated, graded, connected $\F$-algebra. Let $B \subseteq A$ be a graded sub algebra so that $A$ is finitely generated as a left $B$-module. Then $B$ is also a finitely generated algebra. 
\end{lemma}
\begin{proof} Let $\{a_1, \ldots, a_n\}$ be a set of algebra generators of $A$ over $\F$, and let $\{ c_1, \ldots, c_m \}$ be a basis set of $A$ as a left $B$-module, with $c_1 = 1$. We pick a finite set $\{ e_1, \ldots, e_k \} \subseteq B$, so that:
\[ a_j \in \sum_{r,s} \F \, e_r \, c_s, \qquad c_i \, a_j \in \sum_{r,s} \F \, e_r \, c_s, \qquad j \leq n, i \leq m.\] 
Let $b \in B$ be an arbitrary element. Since $A$ is finitely generated, there exists a polynomial $f$ so that $f(a_1,\ldots,a_n) = b$. Using the above properties repeatedly, we may write
\[ b = \sum_i g_i(e_1,\ldots, e_k) \, c_i, \]
for some polynomials $g_i$. But since $A$ is a free left $B$-module, we observe that $g_i=0$ for all $i>1$ and that $b = g_1(e_1,\ldots,e_k)$. It follows that $B$ is generating by $\{e_1, \ldots, e_k\}$. 
\end{proof}
\noindent
As an easy consequence of the above lemmas, we have:
\begin{thm}
Let $\F$ be a field of positive characteristic. Then the dual hopf algebra $A_\F^*$ is a finitely generated $\F$-algebra. It follows that $\Hr_*(\K,\F)$ is also finitely generated $\F$-algebra. 
\end{thm}
\begin{proof} Dualizing the coalgebra structure of $\Hr^*(\KT)$, we observe that $\Hr_*(\KT)$ has the structure of a finitely generated algebra on the set of elements $\delta_{r_i}$ dual to the Schubert basis elements $\delta^{r_i}$. Working with coefficients in a field $\F$ of positive characteristic, the results of previous sections shows that $A_\F^* \subseteq \Hr_*(\KT)$ is a sub algebra. By the first lemma, we see that $\Hr_*(\KT)$ is a finitely generated, free, left $A_\F^*$-module and so the second lemma implies that $A_\F^*$ is a finitely generated algebra. The result about $\Hr_*(\K,\F)$ follows easily once we describe it as an extension of $A_\F^*$ by $\Lambda^*(x_1,\ldots,x_r)$. 
\end{proof}
\section{Examples of rank two}

\noindent
In this section, we describe the structure of the (co)homology of rank two Kac-Moody groups and their flag varieties. By a rank two Kac-Moody group, we shall mean a Kac-Moody group for which the set $I$ has cardinality two. 

\medskip
\noindent Generalized Cartan matrices representing Kac-Moody groups of rank two are given by:
\[
A(a,b) = \left( \begin{matrix} 2 & -a \\
-b & 2 \end{matrix} \right).
\]
\noindent
Throughout this section let $\K=\K(a,b)$ denote the semisimple factor inside the corresponding unitary form. If $ab < 4$, then $\K$ is a compact Lie group. In particular
\begin{alignat*}{2}
(a,b) = (0,0) & \quad \K = SU(2) \times SU(2)\\
(a,b) = (2,1) & \quad \K = \operatorname{Spin}(5) = Sp(2)\\
(a,b) = (1,1) & \quad \K = SU(3) \\
(a,b) = (1,3) & \quad \K = G_2 
\end{alignat*}

\medskip
\noindent
Henceforth, we only work with a generalized Cartan matrix $A = A(a,b)$ with $ab \geq 4$. Let $\To \subset \K$ denote the maximal torus. Then the Weyl group has a presentation given by: $\W = \left< r_1, r_2 |
r_1^2 = r_2^2 = 1\right>$. Thus the Poincar\'{e} series for $\Hr^*(\KT;\Z) $ is
\[
P_t\Hr^*(\KT;{\mathbb Z}) = 1 + 2t^2 + 2t^4 + \ldots = \frac{1+t^2}{1-t^2}
\]
Hence $\Hr^*(\KT,{\mathbb Z})$ contains two elements from the Schubert basis in every positive even degree. Let 
$\delta_n$ be the element $\delta^w$ where $l(w) = n$, $l(wr_1) < l(w)$. Thus $w =\ldots r_1 r_2 r_1$
($n$ terms). Let $\tau_n$ be the other element from the Schubert basis in the same degree. Denote
$\delta_1$ and $\tau_1$ by $\delta$ and $\tau$ respectively. The action of the Weyl group on $\tau$ and $\delta$ can be easily deduced from Claim \ref{claim 2}.

\medskip
\noindent
Given a generalized Cartan metrix $A=A(a,b)$, define integers $c_i$, $d_i$ recursively via:
\[
c_0  = d_0 = 0; \quad c_1 = d_1 =1; \quad c_{j+1} = ad_j - c_{j-1}; \quad d_{j+1} = bc_j - d_{j-1}.
\]

\begin{thm} \label{str} In $\Hr^*(\KT;\Z)$ we have the relations:
\begin{alignat*}{2}
\delta \cup \delta_n &= d_{n+1} \delta_{n+1};  \quad \delta \cup \tau_n = \delta_{n+1} + d_n
\tau_{n+1}\\
\tau \cup \tau_n & = c_{n+1} \tau_{n+1}; \quad \tau \cup \delta_n = \tau_{n+1} + c_n 
\delta_{n+1}
\end{alignat*}
\end{thm}
\begin{proof}
We proceed by induction. Since $\delta_0 = \tau_0 =1$, the result is true for $n=0$. Now write
$\delta \cup \delta_n = A \delta_{n+1} + B \tau_{n+1}$ where $A, B \in \Z$. We recall the annihilation operators $A_1$ and $A_2$ acting diagonally with respect to the Schubert basis. 
Note that $A_1 (\delta \cup \delta_n) = A_1 (A \delta_{n+1} + B \tau_{n+1} )= A \tau_n$, but on the 
other hand we have the twisted derivation property given by theorem \ref{A of a product} and claim \ref{claim 2}:
\[
\begin{split}
 A_1 (\delta \cup \delta_n)& = r_1 (\delta) \cup A_1(\delta_n) + \delta_n \cup A_1 (\delta)\\
& = (\delta - (2 \delta - b \, \tau)) \cup \tau_{n-1} + \delta_n \\
& = b \, \tau \cup \tau_{n-1} - \delta \cup \tau_{n-1} + \delta_n\\
& = bc_n \, \tau_n - \delta_n - d_{n-1} \tau_n + \delta_n\\
& = d_{n+1} \tau_n
\end{split}
\]
where we used induction and the recursive definition of $d_i$.
Thus $A = d_{n+1}$. Now we apply $A_2$ and observe that $B =0$. The other equalities follow 
similarly.
\end{proof}

\begin{defn} Given a generalized cartan matrix $A(a,b)$, we define the generalized binomial coefficients \footnote{This terminology is due to Haynes Miller}:
\begin{align*}
\mbox{D(n,m)} & = \frac{ d_{n+m} d_{n+m -1} \cdots 1}{d_n d_{n-1}\cdots 1 \,\, d_m d_{m-1} \cdots 1}\quad
\text{and}\\ \mbox{C(n.m)} & = \frac{ c_{n+m} c_{n+m -1} \cdots 1}{c_n c_{n-1}\cdots 1 \, \, c_m c_{m-1} \cdots 1}
\end{align*}
Note that if $a = b = 2$, then $c_n = d_n = n$ and thus $\mbox{C(n,m)} = \mbox{D(n,m)} = {n+m \choose n}$.
\end{defn}

\noindent
The previous theorem on the ring structure of $\Hr^*(\KT,\Z)$ immediately implies the following theorem about the cohomology of the partial flag varieties:

\begin{thm} Let $\K_1$, $\K_2$ be the maximal compact subgroups of the standard parabolic subgroups corresponding to $\{1\}$,$\{2\}\subset \{1,2\}$ 
respectively. Then
\begin{alignat*}{2}
\Hr^{2n}(\K/\K_1;\Z) & = \Z \langle \tau_n \rangle;  \quad 
\tau_n \cup \tau_m  =\mbox{C(n,m)} \, \tau_{n+m}\\
\Hr^{2n}(\K/\K_2;\Z)& = \Z \langle \delta_n \rangle;  \quad
\delta_n \cup \delta_m = \mbox{D(n,m)} \, \delta_{n+m}
\end{alignat*}
In particular we see that the generalized binomial coefficients $\mbox{C(n,m)}$ and $\mbox{D(n,m)}$ are integers! \footnote{We thank Kasper Andersen for showing us a nice algebraic proof of integrality.} 
\end{thm}

\begin{remark}
The above theorems completely determine the ring structure of $\Hr^*(\KT; {\mathbb Z})$, and $\Hr^*(\KKI,\Z)$.
\end{remark}
\begin{claim}\label{rational}
$\Hr^*(\KT; {\mathbb Q}) = {\mathbb Q}[\delta, \tau] /J$, where $J$ is the ideal given by the quadratic relation:
$a \delta^2 + b \tau^2 - ab \, \delta \, \tau =0$. In particular $\Hr^*(\KT;{\mathbb Q})$ is generated by $\Hr^2(\KT;{\mathbb Q})$.

\end{claim}
\begin{proof}
In $\Hr^4(\KT; {\mathbb Z})$ we have a relation $a \delta^2 + b \tau^2 - ab \, \delta \, \tau =0$. This yields a map:
\[
{\mathbb Q}[\delta, \tau]/J  \longrightarrow \Hr^*(\KT;{\mathbb Q}).
\]
Using the ring structure of $\Hr^*(\KT,\Q)$, this map is surjective. To see that it is an isomorphism, one simply compares the Poincar\'{e} series. We leave this as an exercise.
\end{proof}
 
\subsection{The additive structure of $\Hr^*(\K,\Z)$, and the Hopf algebras $A_\Z$, $A_{\F_p}$:}
\noindent

\noindent
Now consider the Serre spectral sequence in integral cohomology for the fibration:
\[ \To \longrightarrow \K \longrightarrow \KT. \]
Let $\Hr^1(\To,\Z) = \Z \, \alpha \oplus \Z \, \beta$, with $d_2(\alpha) = \delta$ and $d_2(\beta) = \tau$. 
The ring structure of $\Hr^*(\KT,\Z)$ now allows us to compute the structure of the differential $d_2$. It is easy to see that
$d_2$ is injective on $E_2^{\ast,2}$, and that $E_3 = E_{\infty}$. Let $g_n = \operatorname{gcd}(c_n,d_n)$ denote the g.c.d of the pair $c_n$, $d_n$. The following results are easy consequences of the explicit formulas given in theorem \ref{str}:
\[ E_3^{2n,0} \cong E_3^{2n+2,1} \cong \Z/g_n\Z. \]

\begin{corr} \label{Khom}
The additive structure of $\Hr^*(K,\Z)$ is given by:
\[ \Hr^{2n+3}(\K,\Z) = \Hr^{2n}(\K,\Z) = \Z/g_n \Z. \]
\end{corr}

\noindent
From the above Serre spectral sequence, we notice that $E_3^{\ast,0}$ is given by the Hopf algebra $A_\Z = \Img(\Hr^*(\KT,\Z) \rightarrow \Hr^*(\K,\Z))$, and it can be identified with $\Hr^*(\KT,\Z) \otimes_S \Z$.
Recall that in degree $2n$, $A_\Z$ is a cyclic group of order $g_n$ generated by $\delta_n$ or $\tau_n$:
\[
A_\Z^{2n} = \Z/g_n\Z; \quad g_n = \operatorname{gcd}(c_n,d_n)
\]
The coalgebra structure on $A_\Z$ was induced by: 
\[
\Delta(\delta^w) = \sum \limits_{uv=w} \delta^u \tens \delta^v
\]
where the sum is over all redued expressions, thus
\begin{alignat*}{2}
\Delta(\delta_n) = \sum \limits_{i=0}^n \delta_i \tens \delta_{n-i}'; & \quad
\delta_{n-i}' & = \begin{cases} \delta_{n-i} &\text{if $i$ even}\\
				\tau_{n-i} & \text{if $i$ odd}
		\end{cases}\\
\Delta(\tau_n) = \sum \limits_{i=0}^n \tau_i \tens \tau_{n-i}'; & \quad
\tau_{n-i}' & = \begin{cases} \tau_{n-i} &\text{if $i$ even}\\
				\delta_{n-i} & \text{if $i$ odd}
		\end{cases}
\end{alignat*}
Now fix a prime $p$. Recall that $A_{\F_p} = \F_p \tens_\Z A_\Z$. Hence to understand $A_{\F_p}$, 
we need to know when $p$ divides $g_n$. We have the following theorem on the arithmetic properties of the integers $c_n$ and $d_n$:

\begin{thm}
Let $g_n = \operatorname{gcd}(c_n,d_n)$. Given a prime $p$, there is a smallest positive integer $k$ with the property that $p$ divides $g_k$. Further, 
$p$ divides $g_n$ if and only if $k$ divides $n$. More precisely, $k$ is given by:
\begin{enumerate}
\item[1)] $ k =2p$ if $p$ divides $a$ or $b$ but not both
\item[2)] $ k =p$ if $ab = 4 (\operatorname{mod} p)$ and the conditions of 1) do not hold. In all other cases, we have:
\item[3)] $k = r$ where $r$ is the multiplicative order of any root of the quadratic polynomial given by 
$x^2- (ab-2)x +1$ defined over ${\mathbb F}_{p^2}[x]$.
\end{enumerate}
\end{thm}
\begin{proof}
Note that the $c_i,d_i$ and $g_i$ have terms that look like:

\[
\begin{tabular}{|c|c|c|c|}
\hline
$i$ & $c_i$ & $ d_i $&$ g_i$\\
\hline
$0$ & $0$ & $0 $& $ 0$\\
$1$ & $1$ & $1$&$ 1$\\
$2$ & $a$ & $b$ &$ (a,b)$\\
$3$ & $ab-1$ & $ab-1$&$ ab-1$\\
$4$ & $a(ab-2)$ & $b(ab-2)$&$ (a,b)(ab-2)$\\
\hline
\end{tabular}
\]
\noindent
So the claim can be verified for $p=2$ explicitly. Assuming $p$ is odd, consider the generating 
function
\[
F(x) = \sum \limits_{i=0}^{\infty} {c_i \choose d_i} x^i
\]
This generating function is to be thought of as a formal power series with coefficients in the 
two dimensional vector space over the field ${\mathbb F}_p$. Thus we are interested in when the coefficient
of $x^n$ is zero. Now, we have a functional equation:
\[
\left(x^2 - \left( \begin{matrix}  0 & a \\ b & 0 \end{matrix} \right)x + I \right) F(x) = x 
\left( \begin{matrix} 1\\1 \end{matrix} \right) 
\]
i.e.
\[
(x - M ) (x - M^{-1}) F(x) = x \left( \begin{matrix} 1\\1 \end{matrix}\right)
\]
where $M = \frac{1}{2} \left( \begin{matrix} \mu & a \\b & \mu \end{matrix} \right)$ and 
$\mu = \sqrt{ab-4} \in {\mathbb F}_{p^2}$. 
Thus 

\begin{equation*}
F(x) = \frac{ 1}{(x - M ) ( x - M^{-1})} x \left( \begin{matrix} 1\\1 \end{matrix} \right)
\qquad 
\tag{8}
\end{equation*}

\noindent
First consider the case $ab \equiv 4 (\operatorname{mod}p)$. In this case $M = M^{-1}$, so $(8)$
says
\[
F(x) = \sum \limits_{i=0}^\infty i x^i M^{i+1} \left( \begin{matrix} 1\\1\end{matrix}\right)
\]
Since $M$ is invertible, the coefficients of $x^n$ are zero if and only if $n$ is a multiple
 of $p$, and that is what we wanted to show.
For all other cases $M - M^{-1}$ is invertible, so $(8)$ becomes:
\[
F(x) = \sum \limits_{i=0}^\infty x^i \left( \frac{ M^i - M^{-i}}{M - M^{-1}}\right)
\left( \begin{matrix} 1\\1 \end{matrix} \right)
\]
Thus we are interested in $n$ where 
\[
\left( \frac{ M^n - M^{-n}}{M - M^{-1}}\right)
\left( \begin{matrix} 1\\1 \end{matrix} \right) = \left( \begin{matrix} 0\\0 \end{matrix} \right)
\]
i.e.
\[
M^{2n} \left( \begin{matrix} 1\\1 \end{matrix} \right)= 
\left( \begin{matrix} 1\\1 \end{matrix} \right)
\]
Now consider the case when $p$ divides $a$ or $b$ but not both. Assume without loss of generality 
$p|b$. Then $M = \left( \begin{matrix} \eta & \frac{a}{2} \\ 0 & \eta \end{matrix} \right)$, and
$\eta = \sqrt{-1} \in {\mathbb F}_{p^2}$. Hence 
\[
M^{2n} = \left( \begin{matrix} (-1)^n & (-1)^{n+1} n a \eta \\
0 & (-1)^n\end{matrix} \right)
\]
and $M^{2n} \left( \begin{matrix} 1\\1 \end{matrix} \right) = 
\left( \begin{matrix} 1\\1 \end{matrix} \right)$ if and only if $n$ is a multiple of $2p$.

\noindent In all remaining cases, $M$ is diagonalizable over ${\mathbb F}_{p^2}$. Since $M^2$ has 
determinant $1$, the only way the equality
\[
M^{2n} \left( \begin{matrix} 1\\1 \end{matrix} \right)= 
\left( \begin{matrix} 1\\1 \end{matrix} \right)
\]
can hold is that $M^{2n}$ is the identity matrix. This is equivalent to the condition that the 
eigenvalues of $M^2$ have multiplicative order dividing $n$. These eigenvalues are exactly the 
roots of the characteristic polynomial of $M^2$, which is given by $x^2 - (ab-2)x +1$.
\end{proof}

\medskip
\noindent
From the above theorem, we observe that the Poincar\'{e} series if $A_{\F_p} $ is: 
\[
P_t (A_{\F_p}) = 1 + t^{2k} + t^{4k} + \ldots = \frac{1}{1 - t^{2k}}.
\]
and both $\delta_{mk}$, $\tau_{mk}$ are generators in that degree. It is easier to understand the dual $A_{\F_p}^*$, which turns out to be a polynomial algebra:
\begin{claim}\label{dual}
$A_{\F_p}^*= \F_p[x_{2k}]$ where $x_{2k}$ is a primitive class in degree $2k$.
\end{claim}
\begin{proof}
Let $x \in A_{\F_p}^*$ be any generator in degree $2k$. We prove by induction that $x^n$ is 
a generator in degree $2nk$. Let $\tau_n \in A_{\F_p}$ be some generator in degree $2nk$. We have:
\[
\Delta (\tau_n) = \sum \limits_{i=0}^n \lambda_i \, \tau_i \tens \tau_{n-i}, \quad \lambda_i \neq 0
\]
Now using the induction hypothesis, we get:
\[
\langle x \cdot x^{n-1}, \tau_n \rangle = \langle x \tens x^{n-1}, \Delta \tau_n \rangle = \lambda_1 \, \langle x \tens x^{n-1}, \tau_1 \tens \tau_{n-1} \rangle \neq 0
\]
Thus $x^n$ generates $A_{\F_p}^*$ in degree $2nk$.
\end{proof} 

\begin{thm}
There is an isomorphism of algebras:
\[
\Hr_*(\K; \F_p) = \Lambda(y_3, y_{2k-1}) \tens \F_p[x_{2k}]
\]
where the subscripts denote the homogeneous degree of the generators. These generators are related via a higher Bockstein homomorphism: $\beta^{(m)} x_{2k} = y_{2k-1}$, $m$ being the exponent of $p$ in $g_k$. Moreover, the generators $y_3, y_{2k-1}$ are primitive, and if $p$ is odd, then so is the generator $x_{2k}$.
\end{thm}
\begin{proof}
Recall the short exact sequence of Hopf algebras: 
\[
1 \rightarrow A_{\F_p} \longrightarrow \Hr^*(\K; \F_p) \longrightarrow \Lambda(x_1,x_2) \rightarrow 1
\]
On dualizing, we get:
\[
1 \rightarrow \Lambda(z_1,z_2) \longrightarrow \Hr_*(\K; \F_p) \longrightarrow \F_p[x_{2k}] \rightarrow 1
\]
Now recall that $\Hr_*(\K;\Q) = \Lambda(z)$, $|z| =3$. This forces $|z_1|=3$, where $z_1$ is a permanent 
cycle in the Bockstein spectral sequence for $\Hr_*(\K;\F_p)$ and $|z_2| = 2k-1$, where $z_2$ is the target of a higher Bockstein of height $m$ supported on $x_{2k}$. Let us relabel these classes by their subscript and call them $y_3$ and $y_{2k-1}$ respectively. It is clear for dimensional reasons that these classes are primitive. 
Now since $\F_p[x_{2k}]$ is a free algebra, we may fix a section to the above short exact sequence. Again, for dimensional reasons, the class $x_{2k}$ is primitive, with the possible exception of the case when $p=2$, $k=3$ and the coproduct on $x_{2k}$ is given by: 
\[ \Delta(x_6) = 1\otimes x_6 + x_6 \otimes 1 + x_3 \otimes x_3. \]

\noindent
Finally, to show that $\Hr_*(\K ; \F_p)$ is a tensor product of $\F_p[x_{2k}]$ and $\Lambda(x_1,x_2)$, it is sufficient to show that $[x_{2k}, y_3] = [x_{2k}, y_{2k-1}] = 0$. This is easy to establish since both the elements: $[x_{2k}, y_3]$ and $[x_{2k}, y_{2k-1}]$ are primitive, but on the other hand, there are no non-zero primitive elements in those degrees. 
\end{proof}

\medskip
\noindent
For the sake of completeness, we include the following theorem whose proof can be found in \cite{ACP}. 

\par\medskip

\begin{thm} Given an odd prime $p$, let $\nu_p(s)$ denote the exponent of the prime $p$ dividing $s$. Let $g_k$ denote the first integer in the sequence $\{ g_n\}$, so that $p$ divides $g_k$. Then we have:
\[ \nu_p(g_{sk}) = \nu_p(s) + \nu_p(g_k). \]
\end{thm}

\section{The stable transfer from $\BK$ to $\BN(T)$}
\noindent
It is well known for a compact Lie group $\Gr$, with maximal torus $\To$ and normalizer $\No$, that the suspension spectrum of $\BG_+$ is a stable retract of $\B \No_+$. The retraction is constructed as a transfer map, and uses the essential fact that $\GrN$ is a finite complex (with Euler characteristic equal to one). In the case of Kac-Moody groups $\K$, the space $\KNT$ is not even homologically finite, and so there is no (apriori) obvious transfer map. Nevertheless, in this section we will construct a transfer using the homotopy decomposition of the spaces $\BK$ and $\BN(T)$ given in \cite{Ki,Ki2,BK}. The author would like to thank Bill Dwyer for motivating the argument used in this section. The construction of the stable transfer proceeds as follows:

\medskip
\noindent
Let ${\bf BK}_+$ and ${\bf BN(T)}_+$ denote the suspension spectra of $\BK_+$ and $\BN(T)_+$,  respectively each endowed with a disjoint base point. In order to construct a stable transfer map from ${\bf BK}_+$ to ${\bf BN(T)}_+$, first recall \cite{Ki,Ki2,BK} that the following canonical maps are homotopy equivalences:
\[ \hocolim_{J \in \gC(A)} \BNJ \longrightarrow \BN(T), \qquad \hocolim_{J \in \gC(A)} \BK_J \longrightarrow \BK \]

\noindent
Let us fix a representation $\mathcal{V}$ of the Kac-Moody group $\K$ with a countable basis, and the property that given $J \in \gC(A)$, every representation of $\K_J$ appears in $\mathcal{V}$ with infinite multiplicity. Such a representation is easy to construct: for example, we may take countable sums of all representations of $\K$ of the form ${\mbox L}_{\mu} \otimes {\mbox L}_{\tau}$, where $\mu$ is a dominant weight, $\tau$ is an anti-dominant weight, and ${\mbox L}_{\mu}$ (resp. ${\mbox L}_{\tau}$) denote the highest (resp. lowest) weight irreducible representations of $\K$.

\medskip
\noindent
Let $\mbox {Met}$ denote contravariant functor on the category $\gC(A)$, taking values in spaces given by $\mbox {Met}(J) =  \mbox{Met}_{\K_J} (\mathcal{V})$: the contractible space of $\K_J$-invariant metrics on $\mathcal{V}$. The assignment that sends $J$ to $\mbox{Met}(J)$ is a contravariant functor on $\gC(A)$. Consider the homotopy inverse limit ${\lim}_{\gC(A)}  \mbox{Met}$, of this functor. A standard spectral sequence computes the homotopy groups of ${\lim}_{\gC(A)} \mbox{Met}$ and is given by:
\[ E_2^{i,j} = {\lim}_{\gC(A)}^i \pi_j \mbox{Met}(J) \Rightarrow \pi_{j-i} {\lim}_{\gC(A)} \mbox{Met}. \]
Since $\mbox{Met}(J)$ is contractible for all $J$, we concluse that ${\lim}_{\gC(A)} \mbox{Met}$ is weakly contractible. In particular, it is non-empty. An element in this homotopy inverse limit may be interpreted as a family of metrics parametrized over the simplicial complex $|\gC(A)|$ (given by the geometric nerve of $\gC(A)$) and with the property that the metrics over the face corresponding to the sequence of inclusions: $J_0 < \cdots < J_k$ are $\K_{J_0}$-invariant. Fixing such an element describes $\mathcal{V}$ as a complete universe parametrized over $|\gC(A)|$. Working with this universe throughout, let ${\bf S}^0$ denote the equivariant sphere spectrum. In addition, let ${\bf K_J/N_J(T)}_+$ denote the suspension spectrum of the $\K_J$-space $\KNJT$, endowed with a disjoint base point.

\medskip
\noindent
The subtle part in the construction of the transfer will be to construct a zig-zag of spectra:
\[ \TT: \hocolim_{J \in \gC(A)} {\bf S}^0 \stackrel{\sim}\longleftarrow {\bf X} \longrightarrow \hocolim_{J \in \gC(A)} {\bf K_J/N_J(T)}_+. \]
Each spectrum above will parametrized over the simplicial complex $|\gC(A)|$. The maps in $\TT$ will fiber over self maps of $|\gC(A)|$ which preserves the faces (though not pointwise). By construction, each spectrum above will admit the fiberwise action of the group $\K_{J_0}$ over the face of $|\gC(A)|$ corresponding to the sequence of inclusions $J_0 <\cdots < J_k$. Moreover, $\TT$ will be equivariant with respect to $\K_{J_0}$, over this face. If we let $\TT(J)$ denote the map over the vertex of $|\gC(A)|$ given by the object $J \in \gC(A)$, then our construction will also show that $\TT(J)$ is equivalent to the standard equivariant splitting of ${\bf S}^0$ from ${\bf K_J/N_J(T)}_+$. 

\medskip
\noindent
Taking homotopy orbits of $\TT$, and inverting the equivalence in the zig-zag, we get a map which we will define as the stable transfer $\overline{\TT}: {\bf BK}_+ \rightarrow {\bf BN(T)}_+$:
\[  \overline{\TT}: \hocolim_J {\bf EK}_+ \wedge_{\bf K_J} {\bf S}^0 \longrightarrow \hocolim_J {\bf EK}_+ \wedge_{\bf K_J} ({\bf K_J/N_J(T)}_+). \]

\medskip
\noindent
So it remains to actually construct a map $\TT$ with all the required properties. We begin with some auxiliary constructions. Define ${\bf X}$ as the parametrized spectrum given by the co-end construction induced by obvious restrictions:
\[ {\bf X} = \coprod_{J_0 < \cdots < J_k} \Delta^k \times (\mbox{Emb}(K_{J_k}(A)/N_{J_k}(T))_+ \wedge {\bf S}^0) /\sim ,  \]
where $\mbox{Emb}(\KNJT)$ is the space of $\K_J$-equivariant embeddings of $\KNJT$ in $\mathcal{V}$. The space $\mbox{Emb}(\KNJT)$ is contractible, and so it is clear that the parametrized projection map from ${\bf X}$ to the parametrized spectrum: $\hocolim_{J \in \gC(A)} {\bf S}^0$, is an equivalence. 

\medskip
\noindent
The construction of our stable transfer $\TT$ reduces to the construction of a face-preserving map from ${\bf X}$ to the spectrum $\hocolim_{J \in \gC(A)} {\bf K_J/N_J(T)}_+$, as parametrized family of equivariant splitting maps constructed using the Pontrjagin-Thom collapse construction for equivariant embeddings of $\KNJT$ in $\mathcal{V}$. The subtlety here is to ensure coherence between the individual maps. To address the coherence problem, we will take advantage of the following two general facts (which are easy to prove, and are left to the reader):

\medskip
\begin{enumerate}
\item
Given $ J \in \gC(A)$, let $\gH_J$ denote the Lie algebra: $\gH_J = \langle h \in \gH, \, | \, \alpha_j(h) = 0, \, j \in J \rangle$. Then the infinitesimal action of $\gH_J$ on $\KNST$ is $\K_J$-invariant for any $J < S$. Moreover, the $\gH_J$-fixed set of $\KNST$ is given by $\KNJT$. \\
\item
Fix elements $\rho_J \in \gH_J$ with the property $\alpha_i(\rho_J)=1$ for $i \notin J$. Then the simplicial complex $|\gC(A)|$ given by the geometric nerve of $\gC(A)$ can be canonically identified with an affine subspace of $\gH$, determined by the property that the vertex of $|\gC(A)|$ corresponding to $J \in \gC(A)$, maps to the element $\rho_J$. 
\end{enumerate}

\medskip
\noindent
Now given an embedding $e : \KNJT \rightarrow \mathcal{V}$, recall that the Pontrjagin-Thom construction is given by collapsing a tubular neighbourhood of $e$ to yield an equivariant stable map ${\bf S}^0 \rightarrow {\bf {K_J/N_J(T)}^\eta}$, where ${\bf {K_J/N_J(T)}^\eta}$ is the Thom-spectrum of the stable normal bundle. Let $\tau$ denote the stable tangent bundle of $\KNJT$. Then including $\eta$ into the trivial bundle, we get the $\K_J$-equivariant transfer map given by the composite:
\[ \TT(J,e) : {\bf S}^0 \longrightarrow  {\bf {K_J/N_J(T)}^\eta} \longrightarrow  {\bf K_J/N_J(T)}^{\eta \oplus \tau} = {\bf K_J/N_J(T)}_+. \]
The collection of maps $\TT(J,e)$ yield a map from the restriction of ${\bf X}$ over the zero-skeleton of $|\gC(A)|$ to the corresponding restriction of $\hocolim_{J \in \gC(A)} {\bf K_J/N_J(T)}_+$. It remains to extend this map to the whole simplicial complex $|\gC(A)|$. This is the point where the two properties stated above are crucial. 

\medskip
\noindent
By property (2), we may identify the simplicial complex $|\gC(A)|$ with a piecewise-affine subspace of $\gH$, with the property that the vertex corresponding to $J_i$ is identified with the element $\rho_{J_i}$. Let $\mbox{B}|\gC(A)|$ denote the barycentric subdivision of $|\gC(A)|$. Working inductively with the faces, we may define a simplicial map: $\pi: \mbox{B}|\gC(A)| \rightarrow |\gC(A)|$, with the property $\pi(b(\Delta)) = J_k$, where $b(\Delta)$ denotes the vertex given by the barycenter of a k-dimensional face $\Delta$ of $|\gC(A)|$ corresponding to a sequence of inclusions: $J_0 < \cdots < J_k$. Similarly, define $\lambda$ to be the map $\lambda: \mbox{B}|\gC(A)| \rightarrow |\gC(A)|$ with the property $\lambda(b(\Delta)) = J_0$. 

\medskip
\noindent
Given $\mu \in \Delta$, define $\iota(\mu)$ as the section of the bundle $\eta \oplus \tau \rightarrow \eta$ generated by the vector field given by the infinitesimal action of $\mu \in \gH$. Now consider the following face-preserving map $\TT$ over the k-simplex $\Delta$:
\[ \TT: \Delta \times (\mbox{Emb}(\KNJKT)_+ \wedge {\bf S}^0) \longrightarrow \hocolim_{J \in \gC(A)} {\bf K_J/N_J(T)}_+, \]
\[ 
\TT(x, e) = (\pi(x), \iota(\lambda(x)) \circ \TT(J_k,e)), 
\]
where given $\mu \in \Delta$, the composite $\iota(\mu) \circ \TT(J_k,e)$ is defined as:
\[ {\bf S}^0 \longrightarrow {\bf {K_{J_k}/N_{J_k}(T)}^\eta} \stackrel{\iota(\mu)} \longrightarrow {\bf {K_{J_k}/N_{J_k}(T)}^{\eta \oplus \tau}} = {\bf K_{J_k}/N_{J_k}(T)}_+. \]
Notice that if $J_i$ represents a vertex of $|\gC(A)|$, then property (1) above, implies that the vectorfield on $\KNJKT$, generated by $\rho_{J_i}$ vanishes exactly on $\KNJIT$, and hence the map $\iota(\rho_{J_i})$ collapses the complement of a tubular neighborhood of $\KNJIT$ to the basepoint, thereby showing that the above definition of $\TT$ extends $\TT(J_i,e)$. Moreover, the above definition of $\TT$ is compatible with overlaps of faces, hence $\TT$ yields a face-preserving map of equivariant spectra  from ${\bf X}$ to $\hocolim_{J \in \gC(A)} {\bf K_J/N_J(T)}_+$ over $|\gC(A)|$. The upshot of this argument is that we have a face-preserving transfer map of equivariant spectra over $|\gC(A)|$ given by a zig-zag: 
\[ \TT: \hocolim_{J \in \gC(A)} {\bf S}^0 \stackrel{\sim}\longleftarrow {\bf X} \longrightarrow \hocolim_{J \in \gC(A)} {\bf K_J/N_J(T)}_+. \]

\medskip
\noindent
It is now easy to see that the induced transfer map:
\[  \overline{\TT}: \hocolim_J {\bf EK}_+ \wedge_{\bf K_J} {\bf S}^0 \longrightarrow \hocolim_J {\bf EK}_+ \wedge_{\bf K_J} ({\bf K_J/N_J(T)}_+). \]
is indeed a stable retraction. This can be established by observing that $\overline{\TT}$ induces a map of the respective Bousfield-Kan spectral sequences that compute the stable homotopy of the spectra ${\bf BK}_+$ and ${\bf BN(T)}_+$ respectively. The properties of $\overline{\TT}$ ensure that this map is a retraction on the $E_2$-term, and hence is a retraction. As a consequence, we have:

\begin{thm}
The above map $\overline{\TT}$ is a stable retraction of ${\bf BK}_+$ from ${\bf BN(T)}_+$. In addition, $\overline{\TT}$ is compatible with the stable retractions of each spectra ${\bf BK_J}_+$ from ${\bf BN_J(T)}_+$.
\end{thm}

\begin{remark}
Let ${\mbox Z(K)} \subseteq \To$ denote the center of $\K$. Then the above construction shows that $\overline{\TT}$ is equivariant with respect to the action of $\mbox{BZ(K)}$ on ${\bf BK}_+$ and ${\bf BN(T)}_+$. In particular, $\overline{\TT}$ descends to stable transfers for central quotients of $\K$. 
\end{remark}

\section{Appendix}
\noindent
In this section we establish some basic facts about the topology of the Kac-Moody groups $\K$, and their classifying spaces $\BK$. 

\medskip
\noindent
Recall that a subset $J \subseteq I$ is called spherical if the subgroup $\K_J \subseteq \K$ is a compact Lie group. The poset of spherical subsets of $I$ is denoted by $\gC(A)$. In \cite{KP2} (Theorem A) it is shown that as an abstract group, $\K$ is an amalgamated product of subgroups of the form $\K_J$, where $J \in \gC(A)$ has cardinality at most two. In other words, the following canonical map is an isomorphism:
\[ \colim_{J \in \gC(A), \, |J| \leq 2} \K_J \longrightarrow \K, \]
where the colimit is taken in the category of groups. Now given $J \in \gC(A)$ it is easy to see that $\K_J$ is generated by the groups $\K_j$ for $j \in J$. Hence, the map above factors through a sequence of two surjective maps:
\[ \colim_{J \in \gC(A), \, |J| \leq 2} \K_J \longrightarrow \colim_{J \in \gC(A)} \K_J \stackrel{\varphi} \longrightarrow \K. \]
As a consequence we see that the map $\varphi$ above is an isomorphism of (abstract) groups.

\medskip
\noindent
We now come to the question of topology on $K$. We refer the reader to \cite{K1,Ku} for details regarding the Bruhat decomposition of $\K$ described below. 

\medskip
\noindent
Let $\NoT \subset \K$ denote the normalizer of $\To$. Given $w \in \W$, let $\tilde{w} \in \NoT$ denote any lift of $w$ in $\NoT$. We will denote the space $\B \tilde{w}  \B \cap \K$ by $\Zo_w$. This is a well defined subspace of $\K$ homeomorphic, as a right $\To$-space, to $\C^{l(w)} \times \To$. Now for a generating reflection $r_i$, let $\Y_i \subset \Zo_{r_i}$ be the subspace $\C \times \{ 1 \} \subset \C \times \To$ under the above identification. Then the group product in $\K$ induces a homeomorphism:
\[ \Zo_w = \Y_{i_1} \times \ldots \times \Y_{i_s} \times \To, \]
where $w = r_{i_1} \ldots r_{i_s}$ is a reduced expression. We also have the closure relation:
\[ \overline{\Zo}_w = \coprod_{v \leq w} \Zo_v. \]
With this structure, $\K$ becomes a {\em $\To$-CW complex}, constructed by successively attaching $\To$-cells. The topology is generated by the closed subspaces $\overline{\Zo}_w$. Hence a subspace $\Zo \subseteq \K$ is closed if and only if $\Zo \cap \overline{\Zo}_w$ is closed for all $w \in \W$. 
Now given $J \in \gC(A)$, let $w_0 \in \W_J$ denote the longest element. It follows from the closure relation that $\overline{\Zo}_{w_0} = \K_J$ as compact subspaces of $\K$. 

\medskip
\noindent
Assume now that $H$ is any topological group and that we are given a homomorphism $\phi : \K \rightarrow H$ that restricts to a continuous map on each $\K_J$ for $J \in \gC(A)$. Given an element $w \in \W$, let $w = r_{i_1} \ldots r_{i_s}$ be a reduced expression. Notice that $\phi$ extends to a canonical continuous map $\tilde{\phi}$ from the product $\K_{i_1} \times \ldots \times \K_{i_s}$ to $H$ given by the product of the individual restriction maps. Moreover, $\tilde{\phi}$ factors through the projection map from $\K_{i_1} \times \ldots \times \K_{i_s}$ onto the subspace $\overline{\Zo}_w$. It follows that $\phi$ restricts to a continuous map on $\overline{\Zo}_w$. By the definition of the topology on $\K$, we see and that $\phi$ is in fact a continuous homomorphism. 

\medskip
\noindent
The upshot of the argument given above is that $\K$ is in fact the colimit of the groups $\K_J$ indexed over the poset $\gC(A)$ in the category of {\em Topological Groups}. We conclude:

\begin{thm}
The topological group $\K$ has the following properties:
\begin{enumerate}
\item[(a)] $\K$ is a free $\To$-CW complex of finite type under the right action of $\To$. This structure is compatible with the CW structure on the homogeneous space $\KT$. 
\item[(b)] $\K$ is equivalent to the colimit, in the category of topological groups, of the compact Lie groups $\K_J$ indexed over the poset $\gC(A)$. 
\end{enumerate}
\end{thm}

\begin{remark} Since $\K$ is a $\To$-CW complex, it is built by successively attaching $\To$-cells. Decomposing $\To$ as a CW-complex, we see that $\K$ may be constructed by successively attaching (standard) cells. However, it fails to be a CW complex by virtue of the fact that the boundary of cells being attached may glue to cells of higher dimension. We will call a space built by attaching cells in a possibly non-sequential order a Cell Complex (there is some conflict in the literature on the terminology for such an object). Working inductively with the stages, we see that a cell complex is homotopy equivalent to a CW complex.  
\end{remark}

\begin{remark}
Principal $\K$-bundles that appear in physical applications tend to be defined over base spaces (like solutions of differential equations) that may not have an obviously homotopy type of a CW complex. It is therefore desirable to study the homotopy type of the space ${\mathcal{B}\K}$ that classifies numerable $\K$-bundles. This is Milnor's model \cite{M} of the classifying space of $\K$, which is defined as a colimit of certain spaces $\mathcal{B}_n\K$ under inclusions $\mathcal{B}_n\K \subset \mathcal{B}_{n+1}\K$. The space $\mathcal{B}_n\K$ is defined as the quotient $\mathcal{E}_n\K/\K$, where $\mathcal{E}_n\K$ is the $n$-fold join of $\K$ with itself given the quotient (weak) topology, and the diagonal action of $\K$. Hence $\mathcal{E}_n\K$ can be seen as a quotient of $\Delta^{n-1} \times \K^{\times n}$. Using the fact that $\K$ is a cell complex of finite type, we see that $\mathcal{B}_n\K$ has the structure of a cell complex of finite type. It is clear that the inclusions $\mathcal{B}_n\K \subset \mathcal{B}_{n+1}\K$ are cellular, and therefore, by the previous remark, $\mathcal{B}\K$ has the homotopy type of a CW complex. 
\end{remark}

\end{document}